\documentclass[a4paper,10pt, final]{article}

\usepackage{amsthm,amsfonts,amsmath}
\usepackage{graphicx,psfrag}
\usepackage{hyperref}
\usepackage{a4wide}
\usepackage{enumerate}

\newtheorem{proposition}{Proposition}[section]
\newtheorem{theorem}[proposition]{Theorem}

\newtheorem{definition}[proposition]{Definition}

\newenvironment{proofof}[1]{\smallskip\noindent{\textbf{Proof~of~#1.}}%
  \hspace{1pt}}{\hspace{-5pt}{\nobreak\quad\nobreak\hfill\nobreak%
    $\square$\vspace{2pt}\par}\smallskip\goodbreak}

\numberwithin{equation}{section}
\allowdisplaybreaks[4]

\renewcommand{\phi}{\varphi}
\renewcommand{\epsilon}{\varepsilon}
\renewcommand{\theta}{\vartheta}
\renewcommand{\L}[1]{\mathbf{L^#1}}
\newcommand{\Lloc}[1]{\mathbf{L^{#1}_{loc}}}
\newcommand{\C}[1]{\mathbf{C^{#1}}}
\newcommand{\Cc}[1]{\mathbf{C_c^{#1}}}

\newcommand{\W}[2]{\mathbf{W^{#1,#2}}}
\newcommand{\BV}{\mathbf{BV}}
\newcommand{\modulo}[1]{{\left|#1\right|}}
\newcommand{\norma}[1]{{\left\|#1\right\|}}
\newcommand{\reali}{{\mathbb{R}}}

\newcommand{\tv}{\mathop\mathrm{TV}}
\renewcommand{\O}{\mathinner{\mathcal{O}(1)}}

\newcommand{\spt}{\mathop\mathrm{spt}}

\renewcommand{\d}[1]{\mathinner{\mathrm{d}{#1}}}
\renewcommand{\div}{\mathop{\rm{div}}\nolimits_x} %{\nabla\cdot}%
\renewcommand{\hat}[1]{\widehat{#1}}
\newcommand{\gradx}{\mathop{\rm{grad}}\nolimits_x}
\newcommand{\gradA}{\mathop{\rm{grad}}\nolimits_A}

\newcommand{\Caption}[1]{
  \begin{minipage}{0.75\linewidth}
    \caption{\small{#1}}
  \end{minipage}}

\begin{document}

\title{NonLocal Systems of Balance Laws \\ in Several Space Dimensions
  \\ with Applications to Laser Technology}

\author{Rinaldo M. Colombo$^1$ \and Francesca Marcellini$^2$}

\footnotetext[1]{INDAM Unit, University of Brescia,
  Italy. \texttt{rinaldo.colombo@unibs.it}}

\footnotetext[2]{Dept.~of Mathematics and Applications, University of
  Milano-Bicocca, Italy. \texttt{francesca.marcellini@unimib.it}}

\maketitle

\begin{abstract}
  \noindent For a class of systems of nonlinear and nonlocal balance
  laws in several space dimensions, we prove the local in time
  existence of solutions and their continuous dependence on the
  initial datum. The choice of this class is motivated by a new model
  devoted to the description of a metal plate being cut by a laser
  beam. Using realistic parameters, solutions to this model obtained
  through numerical integrations meet qualitative properties of real
  cuts.  Moreover, the class of equations considered comprises a model
  describing the dynamics of solid particles along a conveyor belt.

  \medskip


  \medskip

\end{abstract}

\section{Introduction}
\label{sec:Intro}

We are concerned with a system of $n$ balance laws in several space
dimensions of the type
\begin{equation}
  \label{eq:2}
  \left\{
    \begin{array}{l}
      \partial_t u_i
      +
      \div  \phi_i (t,x,u_i, \theta * u)
      = \Phi_i (t,x,u_i, \theta*u)
      \\
      u_i (0,x) = \bar u_i (x)
    \end{array}
  \right.
  \qquad i= 1, \ldots, n \,.
\end{equation}
Here, $t \in \left[0, +\infty\right[$ is time, $x \in \reali^N$ is the
space coordinate and $u \equiv (u_1, \ldots, u_n)$, with $u_i = u_i
(t,x)$, is the unknown. The function $\theta$ is a smooth function
defined in $\reali^N$ attaining as values $m \times n$ matrices, so
that
\begin{displaymath}
  \theta \in \Cc2 (\reali^N; \reali^{m\times n})
  \,,\qquad
  \left(\theta * u (t)\right) (x)
  =
  \int_{\reali^N} \theta (x-\xi) \; u (t, \xi) \d\xi
  \,,\qquad
  \left(\theta * u (t)\right) (x) \in \reali^m \,.
\end{displaymath}
The flow $\phi \equiv (\phi_1, \ldots, \phi_n)$, with $\phi_i
(t,x,u_i, A) \in \reali^N$, and the source $\Phi \equiv (\Phi_1,
\ldots, \Phi_n)$, with $\Phi_i (t,x,u_i, A) \in \reali$, have the
peculiar property that the equations are coupled only through the
nonlocal convolution term $\theta * u$.

The driving example for our considering the class~\eqref{eq:2} is a
new model for the cutting of metal plates by means of a laser beam,
presented in Section~\ref{sec:M}. A sort of \emph{pattern formation}
phenomenon, typical of various nonlocal
equations~\cite{ColomboGaravelloLecureux2012}, accounts for the
formation of the well known \emph{ripples} whose insurgence deeply
affects the quality of the cuts. In fact, two type of lasers are
mainly used in the cutting of metals: $CO_2$ lasers and fiber
lasers. The former ones are more powerful and more precise, but also
more expensive. Recent technological improvements are apparently going
to allow also to the cheaper devices of the latter type to cut thick
plates, nowadays treated typically with $CO_2$ lasers. Unfortunately,
a typical drawback of fiber lasers is that along the cut ripples are
generated, see~\cite{Schulz1997, schultz, VossenSchuttler}. The
modeling of these ripples often relies on the introduction of
\emph{imperfections} in the metal or of \emph{inaccuracies} in the
laser management, see also~\cite{0022-3727-32-11-307,
  Schulz2009}. Here, using realistic numeric parameters, we obtain the
formation of a geometry similar to the ripples observed in industrial
cuts. We remark that in the present construction neither the initial
data nor the parameters in the equations contain any oscillating term.

Furthermore, in Section~\ref{sec:CB}, we slightly extend the model
introduced in~\cite{SchleperGoettlich} to describe the dynamics of
bolts along a conveyor belt. The resulting equations fit in the
present framework and is proved to be well posed.

Besides, we also note that several crowd dynamics models considered in
the literature fit into~\eqref{eq:2},
e.g.~\cite{ColomboGaravelloLecureux2012, ColomboLecureuxPerDafermos,
  DiFrancescoEtAl, Hughes2002}.

The particular structure of~\eqref{eq:2} allows to prove its well
posedness. Indeed, for small times, system~\eqref{eq:2} admits a
unique solution $u = u (t,x)$. Moreover, $u$ is proved to be a
continuous function of time with respect to the $\L1$ topology and an
$\L1$--Lipschitz continuous function of the initial datum $\bar u$. In
all this, the particular coupling among the equations in~\eqref{eq:2}
plays a key role. At present, the well posedness of general systems of
balance laws in several space dimensions is a formidable open
problem. In the present work, the functional setting is provided by
$\L1 \cap \L\infty \cap \BV$, as usual in the framework of nonlocal
conservation laws. The existence result is obtained through a careful
use of the general estimates~\cite{ColomboMercierRosini,
  MagaliV2}. They provide the necessary analytic tool to apply Banach
Contraction Theorem.

A preliminary result related to Theorem~\ref{thm:main} below is
presented for instance in~\cite{AggarwalColomboGoatin}, see
also~\cite{KarlsenNonLocal}. There, the existence of solution
to~\eqref{eq:2} in the case $\Phi \equiv 0$ is obtained proving the
convergence (up to a subsequence) of a Lax--Friedrichs type
approximate solutions. Note however that differently from the present
situation, in the case considered in~\cite{AggarwalColomboGoatin},
positive initial data yield positive solutions so that the $\L1$ norm
is conserved.

We remark that most of the results related to nonlocal balance laws
are currently devoted to conservation laws, i.e., to equations that
lack any source term. Here, we allow for the presence of source terms
that can be nonlinear in both the unknown variable $u$ and the
convolution term $\theta * u$. The unavoidable cost of this extension
is a local in time existence result, as shown by an example in
Section~\ref{sec:AR}.

Nonlocal conservation and balance laws are currently widely considered
in various modeling frameworks. Besides those of crowd dynamics, laser
cutting and conveyor belt dynamics considered above, we recall for
instance granular materials, see~\cite{AmadoriShen}, and vehicular
traffic, see~\cite{BlandinGoatin}. For a different approach, based on
measure valued balance laws, we refer to~\cite{PiccoliTosinARMA}.

\smallskip

The paper is organized as follows: the next section is devoted to the
analytic results. Section~\ref{sec:M} presents the laser cutting
model, its well posedness and some qualitative properties with the
help of numerical integrations. Conveyor belts dynamics is the subject
of Section~\ref{sec:CB}. All analytic proofs are postponed to the last
Section~\ref{sec:TD}.

\section{Analytic Results}
\label{sec:AR}

Throughout, we denote by $\gradx f$, respectively $\div f$, the
gradient, respectively the divergence, of $f$ with respect to the $x$
variable, with $x \in \reali^N$. All norms in function spaces are
denoted with a subscript indicating the space, as for instance in
$\norma{u (t)}_{\L1 (\reali^N; \reali^n)}$. When no space is
indicated, the norm is the usual Euclidean norm in $\reali^k$, for a
suitable $k$, as for instance in $\norma{u (t,x)}$. Throughout, we fix
the non trivial time interval $ \hat I = [0, \hat T]$. For any $U >
0$, we also denote $\mathcal{U}_U = [-U, \, U]$.

Our starting point is the definition of solution to~\eqref{eq:2},
which extends~\cite[Definition~2.1]{AggarwalColomboGoatin} to the case
of balance laws.

\begin{definition}
  \label{def:sol}
  Fix a positive $T$. Let $\bar u\in \L\infty (\reali^N,\reali^n)$. A
  map $u:[0,T]\rightarrow \L\infty (\reali^N,\reali^n)$ is a solution
  on $[0,T]$ to~\eqref{eq:2} with initial datum $\bar u$ if, for $i=
  1, \ldots, n$, setting for all $w\in \reali$
  \begin{displaymath}
    \widetilde{\phi_i} (t,x,w)
    =
    \phi_i \left(t,x,w,(\theta * u)(t,x)\right)
    \qquad \mbox{ and } \qquad
    \widetilde\Phi_i(t,x,w)
    =
    \Phi_i\left(t,x,w,(\theta*u)(t,x) \right)
  \end{displaymath}
  the map $u$ is a Kru\v zkov solution to the system
  \begin{equation}
    \label{eq:sol}
    \left\{
      \begin{array}{l}
        \partial_t u_i
        +
        \div  \widetilde\phi_i (t,x,u_i)
        = \widetilde\Phi_i(t,x,u_i)
        \\
        u_i (0,x) = \bar u_i (x)\,
      \end{array}
    \right.
    \qquad i= 1, \ldots, n \,.
  \end{equation}
\end{definition}

\noindent Above, for the definition of Kru\v zkov solution we refer to
the original~\cite[Definition~1]{Kruzkov}.

We are now ready to state the main result of the present paper.

\begin{theorem}
  \label{thm:main}
  Assume that there exists a function $\lambda \in (\C0 \cap \L1)
  (\hat I \times \reali^N \times \reali^+; \reali^+)$ such that:
  \begin{description}
  \item[($\boldsymbol{\phi}$)] For any $U > 0$, $\phi \in (\C2 \cap
    \W{2}{\infty}) (\hat I \times \reali^N \times \mathcal{U}_U \times
    \mathcal{U}_U^m; \reali^{n\times N})$ and for all $t \in \hat I$,
    $x \in \reali^N$, $u \in \mathcal{U}_U$, $A \in \mathcal{U}_U^m$
    \begin{displaymath}
      \max
      \left\{
        \begin{array}{ll}
          \norma{\gradx \phi (t,x,u,A)} \,,
          &
          \norma{\div \phi (t,x,u,A)} \,,
          \\[3pt]
          \norma{\gradx \div \phi (t,x,u,A)} \,,
          &
          \norma{\gradx \gradA \phi (t,x,u,A)} \,,
          \\[3pt]
          \norma{\gradA \phi (t,x,u,A)} \,,
          &
          \norma{\gradA^2 \phi (t,x,u,A)}
        \end{array}
      \right\}
      \leq
      \lambda (t,x,U) \,.
    \end{displaymath}
  \item[($\boldsymbol{\Phi}$)] For any $U > 0$, $\Phi \in (\C1 \cap
    \W{1}{\infty}) (\hat I \times \reali^N \times \mathcal{U}_U \times
    \mathcal{U}_U^m; \reali^n)$ and for all $t \in \hat I$, $x \in
    \reali^N$, $u \in \mathcal{U}_U$, $A \in \mathcal{U}_U^m$
    \begin{displaymath}
      \max
      \left\{
        \norma{\Phi (t,x,u,A)}
        \,,\;
        \norma{\gradx \Phi (t,x,u,A)}
      \right\}
      \leq
      \lambda (t,x,U) \,.
    \end{displaymath}
  \item[($\boldsymbol{\theta}$)] $\theta \in \Cc2 (\reali^N;
    \reali^{m\times n})$.
    % \item[($\boldsymbol{\bar u}$)] $\bar u \in (\L1 \cap \L\infty
    %   \cap
    %   \BV) (\reali^N; \reali^n)$.
  \end{description}
  \noindent Then, for any positive $\bar{\mathcal{C}}$ there exists a
  positive $T_* \in I$ and positive $\mathcal{L}$, $\mathcal{C}$ such
  that for any datum
  \begin{equation}
    \label{eq:barU}
    \bar u \in (\L1 \cap \L\infty \cap \BV) (\reali^N; \reali^n)
    \mbox{ with }
    \norma{\bar u_i}_{\L1 (\reali^N; \reali^n)} \leq \bar{\mathcal{C}} ,\;
    \norma{\bar u_i}_{\L\infty (\reali^N; \reali^n)} \leq \bar{\mathcal{C}}
    \mbox{ and }
    \tv (\bar u_i) \leq \bar{\mathcal{C}},
  \end{equation}
  problem~\eqref{eq:2} admits a unique solution
  \begin{displaymath}
    u \in \C0 \left([0,T_*]; \L1 (\reali^N; \reali^n)\right)
  \end{displaymath}
  in the sense of Definition~\ref{def:sol}, satisfying the bounds
  \begin{displaymath}
    \norma{u (t)}_{\L1 (\reali^N; \reali^n)} \leq \mathcal{C} \,,\;\;
    \norma{u (t)}_{\L\infty (\reali^N; \reali^n)} \leq \mathcal{C}
    \mbox{ and }
    \tv (u (t)) \leq \mathcal{C}\,,
  \end{displaymath}
  for all $t \in [0,T_*]$. Moreover, if also $\bar w$
  satisfies~\eqref{eq:barU} and $w$ is the corresponding solution
  to~\eqref{eq:2}, the following Lipschitz estimate holds:
  \begin{displaymath}
    \norma{u (t) - w (t)}_{\L1 (\reali^N; \reali^n)}
    \leq
    \mathcal{L} \;
    \norma{\bar u - \bar w}_{\L1 (\reali^N; \reali^n)} \,.
  \end{displaymath}
\end{theorem}

\medskip

\noindent The proof is deferred to Section~\ref{sec:TD}. Observe that
the whole construction in the present paper can be easily extended
substituting the convolution $\theta * u$ with a nonlocal operator
having suitable properties that comprise those of the convolution, as
was done for instance in~\cite{ColomboGaravelloLecureux2012,
  ColomboLecureuxPerDafermos}.

\medskip

\begin{figure}[!h]
  \centering
  \includegraphics[width=0.4\textwidth,trim=45 20 30 45]{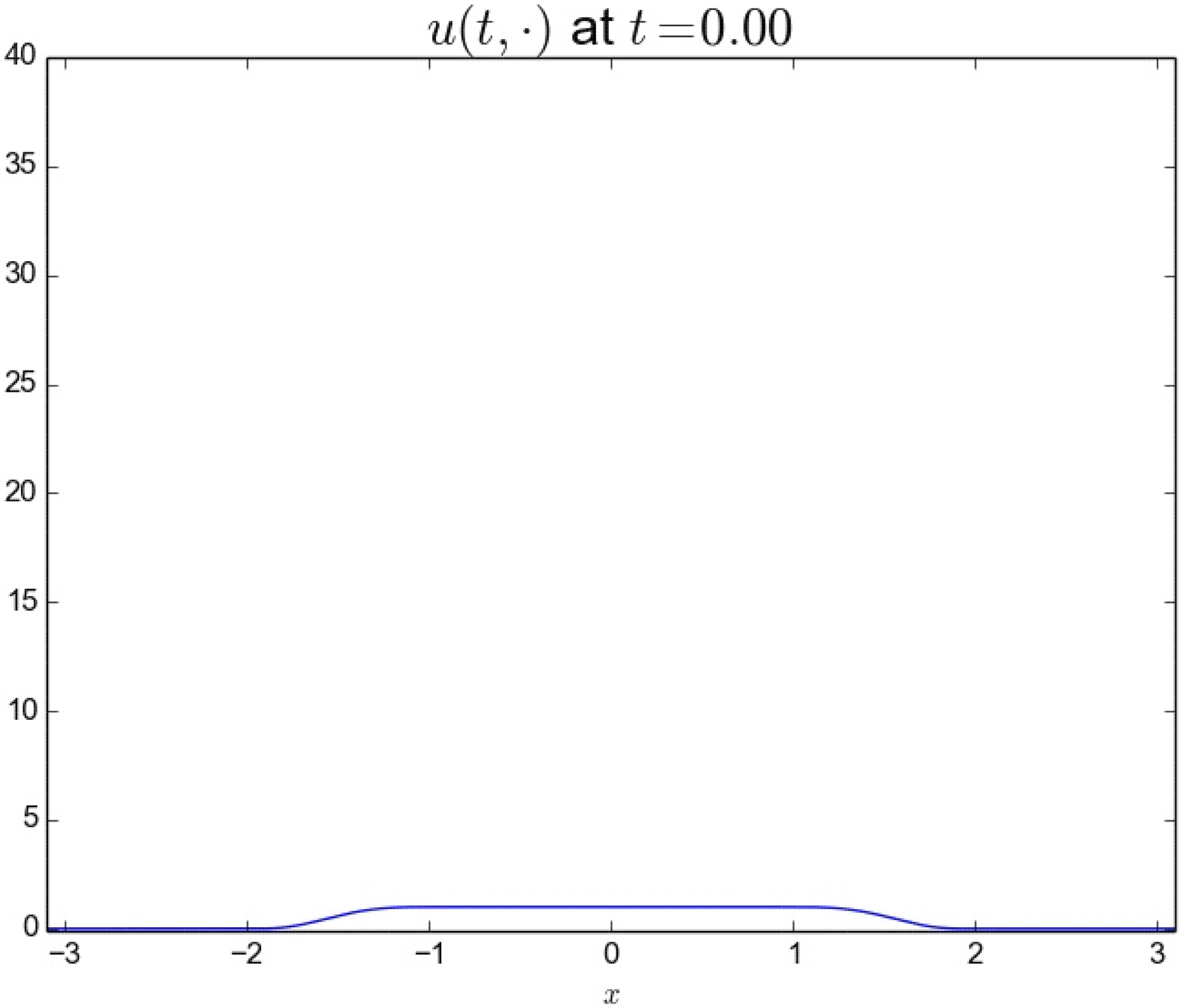}%
  \includegraphics[width=0.4\textwidth,trim=45 20 30 45]{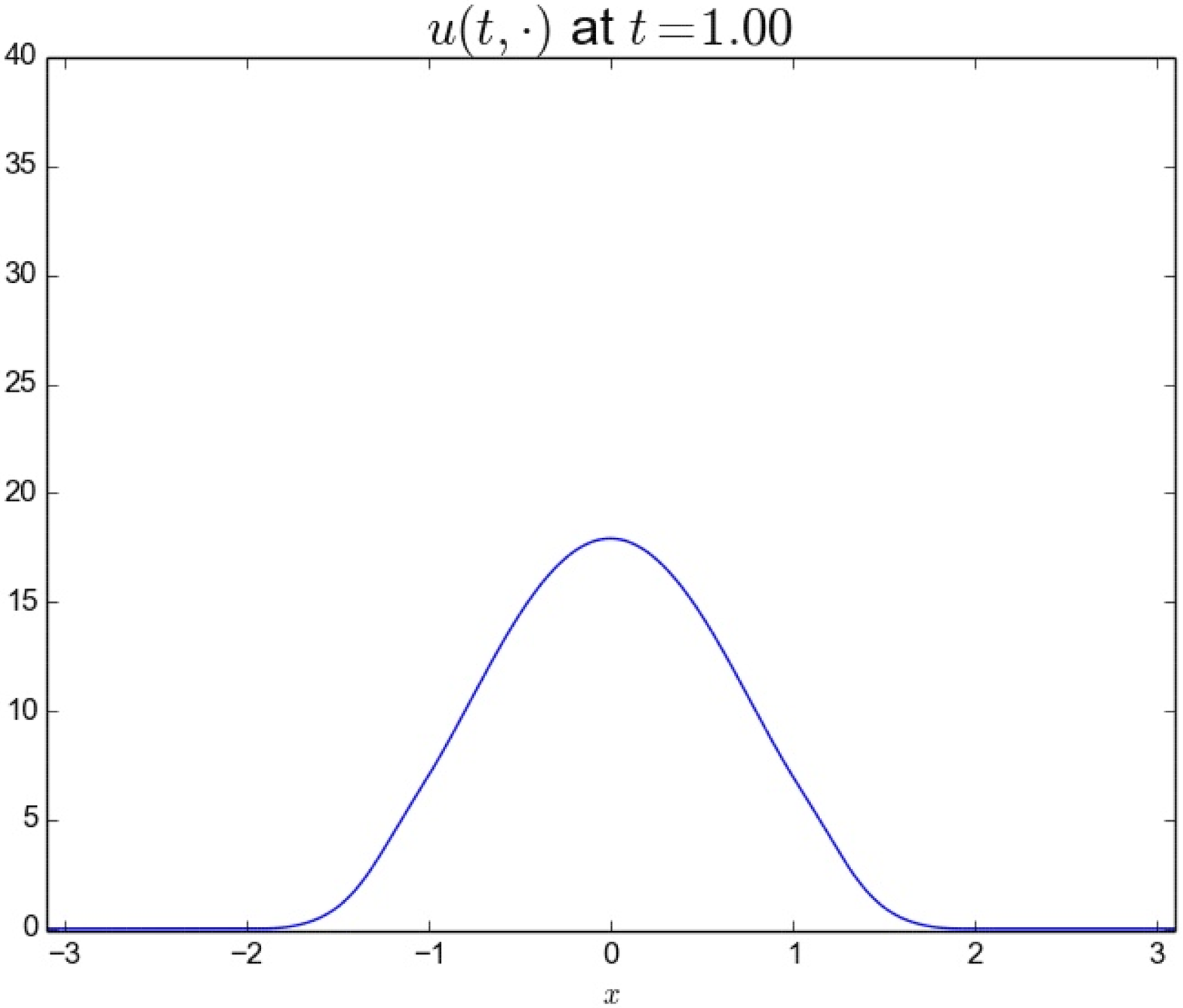}\\
  \includegraphics[width=0.4\textwidth,trim=45 20 30 45]{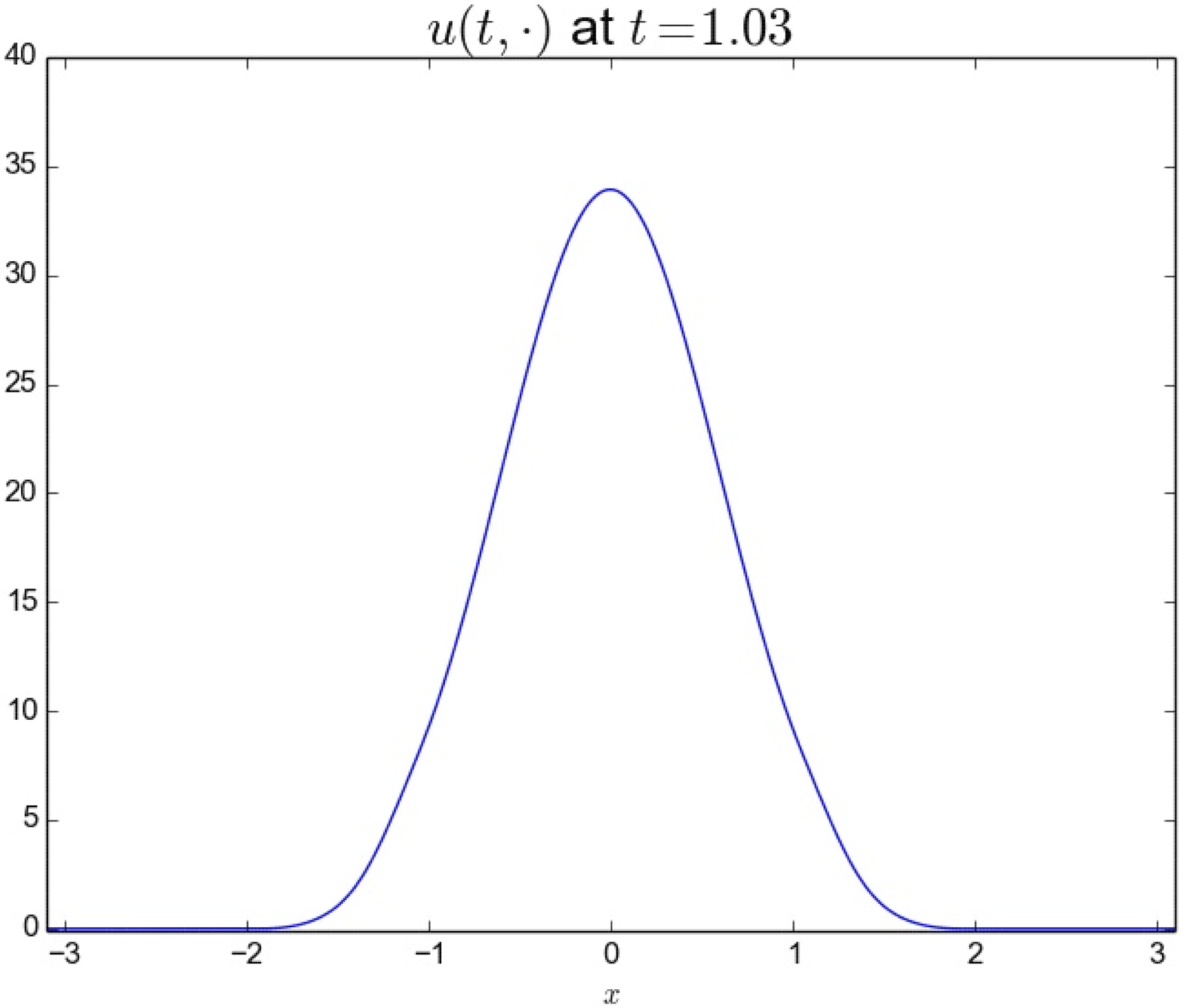}%
  \includegraphics[width=0.4\textwidth,trim=45 20 30 45]{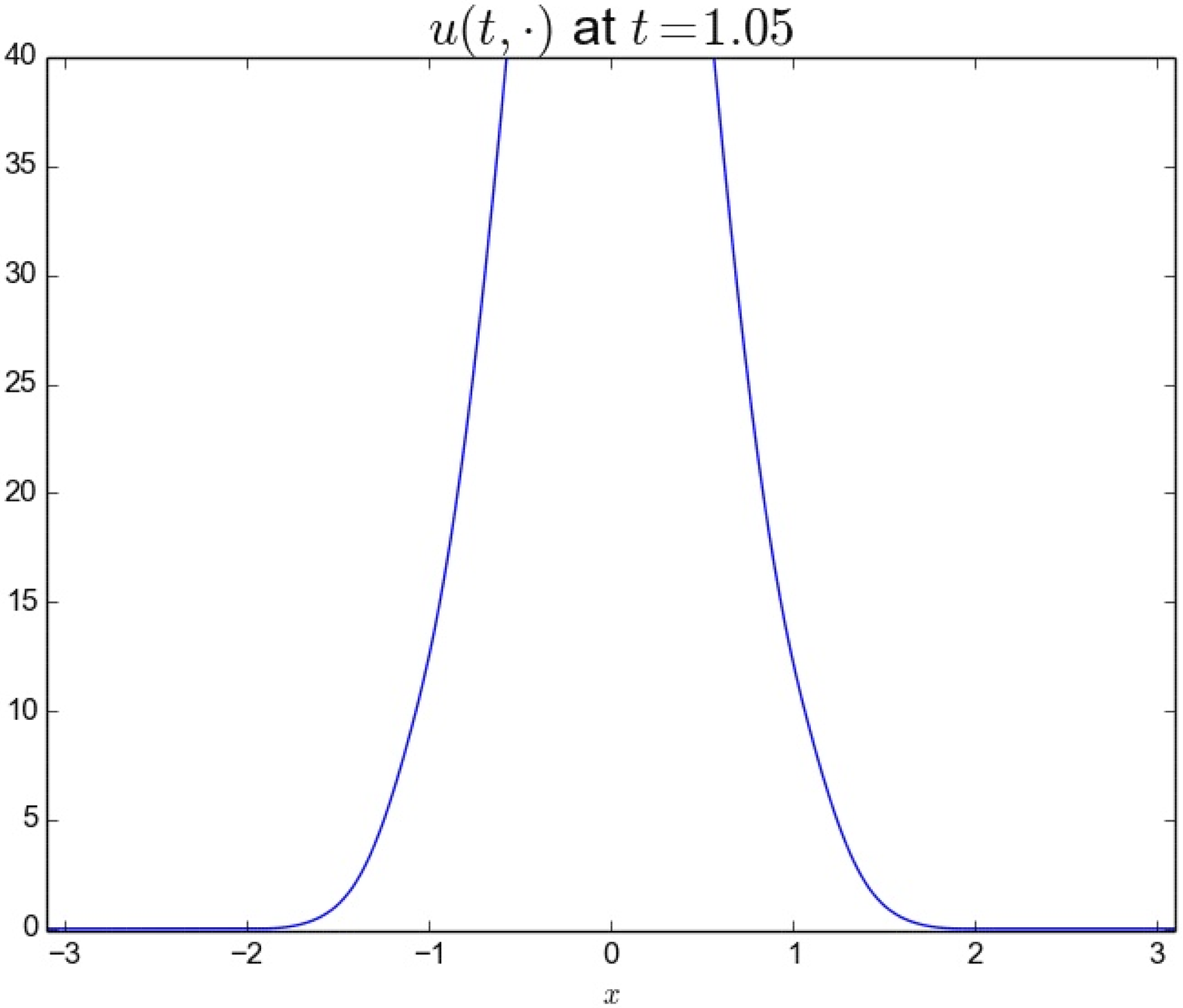}%
  \Caption{Numerical integration of~\eqref{eq:eseSi} for $t \in [0,\,
    1.05]$. The values of the $\L1$ norm of the solution is plotted
    vs.~time in Figure~\ref{fig:eseiu}.}
  \label{fig:ese}
\end{figure}

\begin{figure}[!h]
  \centering
  \includegraphics[width=0.5\textwidth]{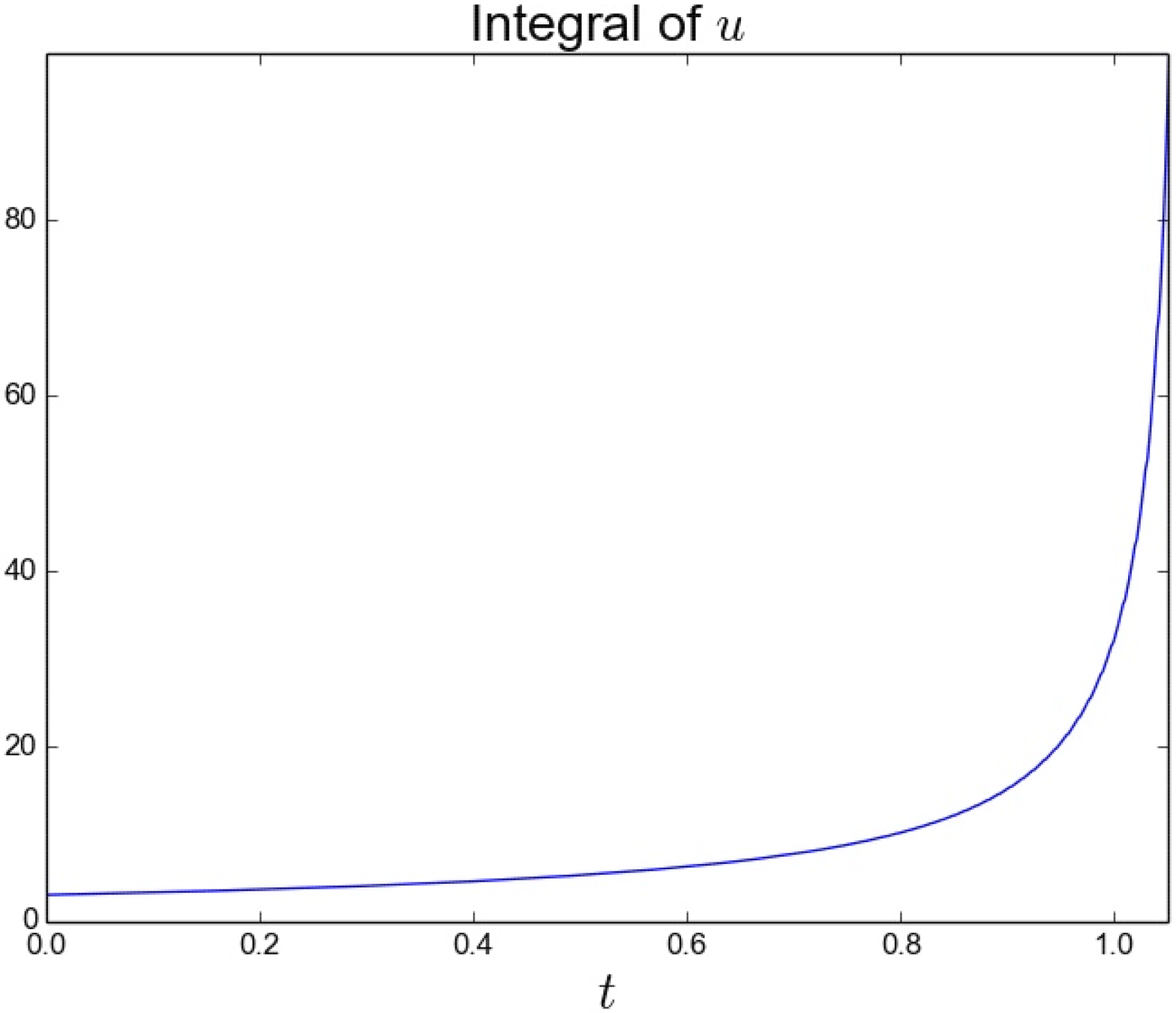}
  \Caption{$\L1$ norm of the solution to~\eqref{eq:eseSi}, suggesting
    a blow up at finite time, similar to the solution
    to~\eqref{eq:eseNo}.}
  \label{fig:eseiu}
\end{figure}

A natural question arises, namely whether the above result can be
extended to ensure the global in time existence of solutions. In this
connection, consider the following particular case of~\eqref{eq:2}
\begin{equation}
  \label{eq:eseNo}
  \left\{
    \begin{array}{l}
      \partial_t u = (u * \eta) \, u
      \\
      u (0,x) = 1 \,.
    \end{array}
  \right.
\end{equation}
Here, $n = 1$ and $m=1$ while $N$ does not play any particular
role. Moreover, $\eta \in \Cc2 (\reali^N; \reali)$ is non negative and
satisfies $\int_{\reali^N} \eta (x) \d{x} = 1$. The solution is $u
(t,x) = 1/ (1-t)$, which exists only up to time $t = 1$. The above
example~\eqref{eq:eseNo} admits an explicit solution but does not fit
into the setting of Theorem~\ref{thm:main}, since the initial datum is
not in $\L1(\reali^N; \reali)$. On the other hand, setting $N=1$, the
similar problem
\begin{equation}
  \label{eq:eseSi}
  \begin{array}[c]{c}
    \left\{
      \begin{array}{l}
        \partial_t u = (u * \eta) \, u \, \psi (x)
        \\
        u (0,x) = \psi (x)
      \end{array}
    \right.
    \quad \mbox{ where }
    \\[30pt]
    \psi (x)
    =
    \left\{
      \begin{array}{lrcl}
        1 & \modulo{x} & \in &[0, 1]
        \\
        \left(1- (x-1)^3\right)^4 & \modulo{x} & \in & \left]1, \, 2\right[
        \\
        0 & \modulo{x}& \in & \left[2, +\infty\right[
      \end{array}
    \right.
  \end{array}
  \begin{array}{c}
    \includegraphics[width=0.4\textwidth]{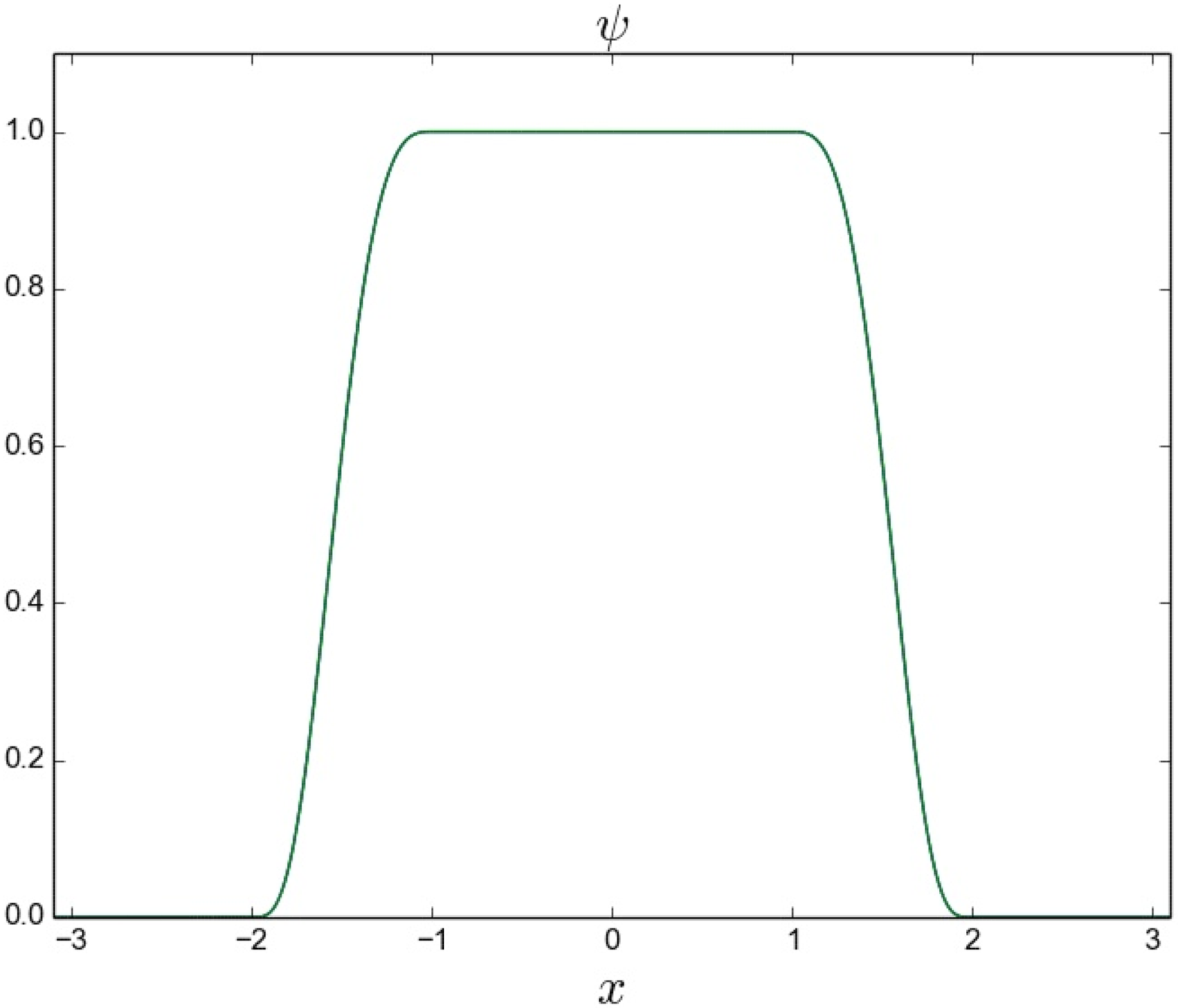}
  \end{array}
\end{equation}
apparently has a qualitatively analogous blow up pattern, as shown by
the numerical integration displayed in Figure~\ref{fig:ese}. To obtain
it, we use an explicit forward Euler method, with space mesh $\Delta x
= 10^{-3}$ and time mesh $\Delta t = 10^{-3}$ on the space domain
$[-3, \, 3]$ and for $t \in [0, \, 1.05]$.
The graph of the $\L1$ norm of the numerical solution
to~\eqref{eq:eseSi} is in Figure~\ref{fig:eseiu}.
It is straightforward to see that~\eqref{eq:eseSi} fits into the
framework of Theorem~\ref{thm:main}, setting
\begin{displaymath}
  \begin{array}{rcl}
    N & = & 1
    \\
    n & = & 1
    \\
    m & = & 1
  \end{array}
  \qquad\qquad
  \begin{array}{rcl}
    \phi (t,x,u,A) & = & 0
    \\
    \Phi (t,x,u,A) & = & \psi (x) \, u \, A \,.
  \end{array}
\end{displaymath}
The requirements~\textbf{($\boldsymbol{\phi}$)}
and~\textbf{($\boldsymbol{\Phi}$)} are easily seen to be satisfied.

\section{A Laser Beam Cutting a Metal Plate}
\label{sec:M}

\begin{figure}[!h]
  \centering
  \begin{psfrags}
    \psfrag{x}{$x_1$} \psfrag{y}{$x_2$} \psfrag{z}{$z$}
    \includegraphics[width=5cm]{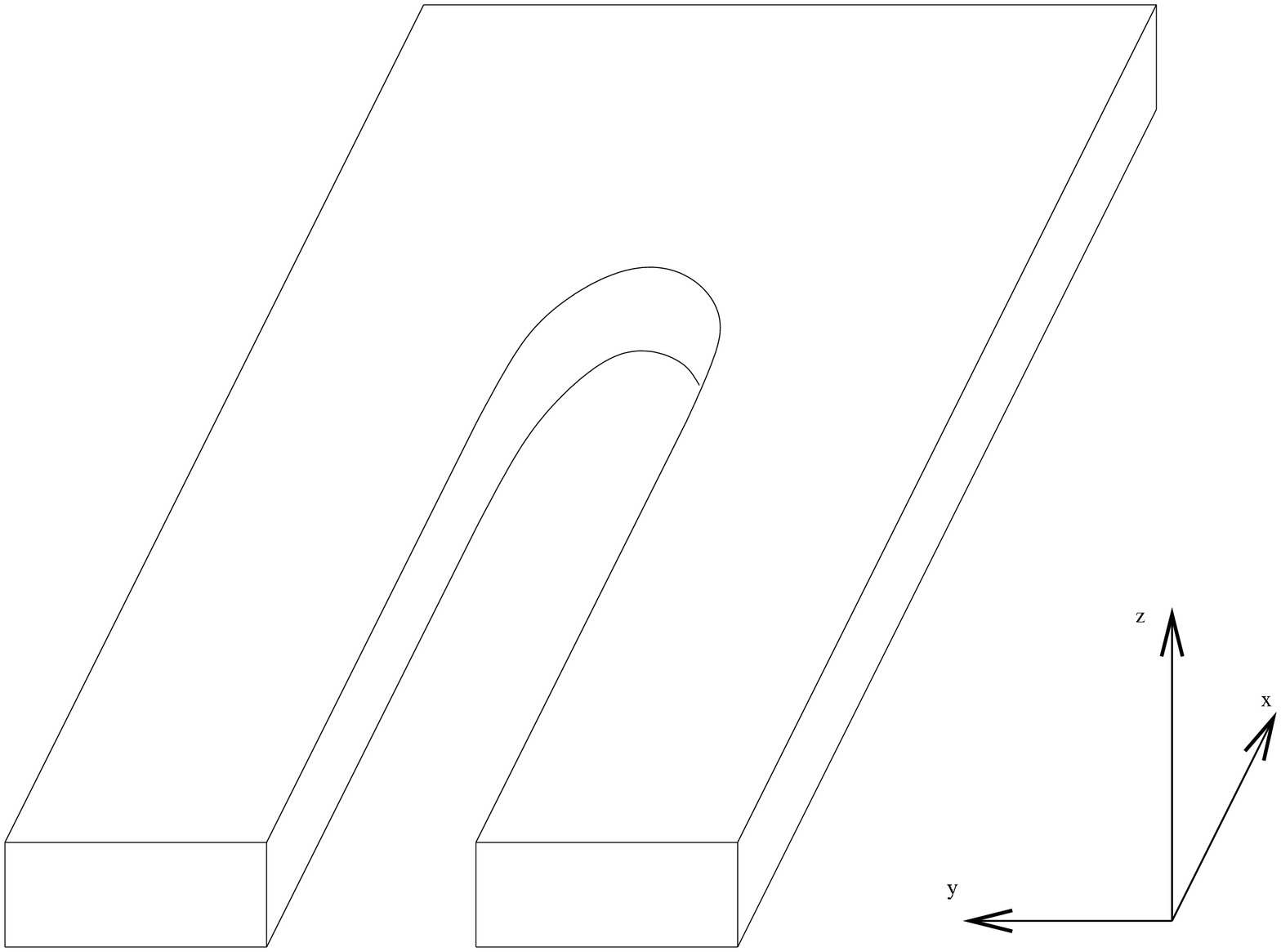}
  \end{psfrags}
  \hfil\hfil
  \begin{psfrags}
    \psfrag{x}{$x_1$} \psfrag{z}{$z$} \psfrag{hs}{$h_s$}
    \psfrag{hm}{$h_m$}
    \includegraphics[width=5cm]{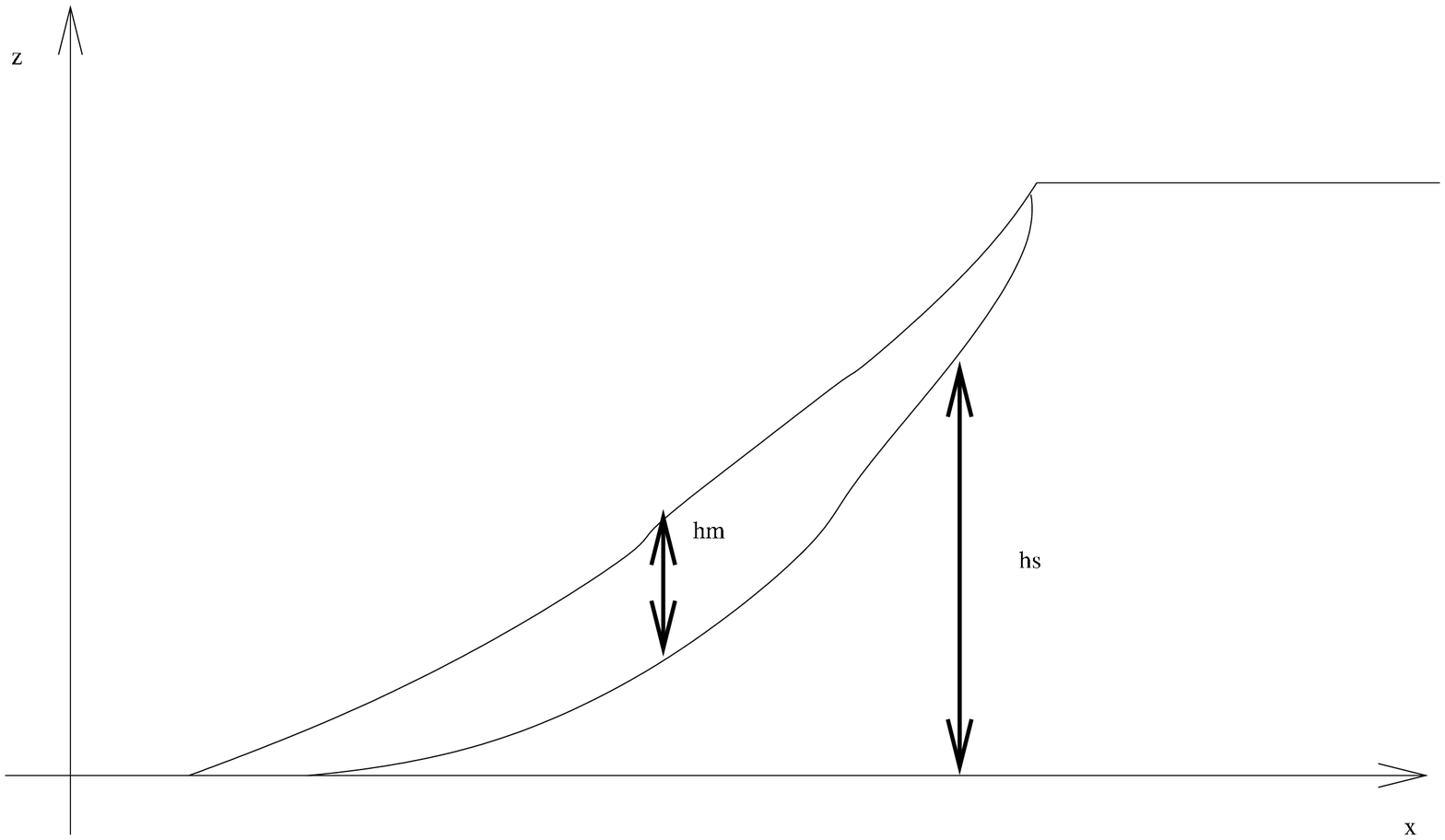}
  \end{psfrags}
  \Caption{Left, reference frame with respect to the metal plate being
    cut. The laser beam is parallel to the $z$ axis, while the plate
    lies on the $z=0$ plain. Right, the distinction between the melted
    part $h_m$ and the solid one $h_s$.\label{fig:SteelGeo}}
\end{figure}

A thin horizontal metal plate can be cut by means of a moving vertical
laser beam. More precisely, the laser energy melts the metal along a
prescribed trajectory. A wind, suitably provoked around the beam,
pushes the melted material downwards. For its industrial interest,
this phenomenon is widely considered in the specialized literature,
see~\cite{Gross2006, Gross2003, Gross2005, 0022-3727-44-10-105502,
  Schulz1997, Schulz2009, schultz, 0022-3727-20-9-016,
  VossenSchuttler}, while information specific to the cut of aluminum
are for instance in~\cite{Smith}. A phenomenological description of
the whole process can be summarized as follows. We fix a 3D geometric
framework, with the laser beam parallel to the vertical $z$ axis, see
Figure~\ref{fig:SteelGeo}, left. The trajectory of the laser is
prescribed by the map $x_L = x_L (t)$.
We distinguish the height $h_s$ of the solid metal and that of the
melted part, denoted $h_m$, see Figure~\ref{fig:SteelGeo}, right.

\begin{figure}[!h]
  \centering
  \begin{psfrags}
    \psfrag{x}{$x_1$} \psfrag{z}{$z$} \psfrag{H}{$z = H (t,x)$}
    \psfrag{s}{$z = h_s (t,x)$} \psfrag{l}{Laser axis}
    \psfrag{a}{$\alpha$} \psfrag{hm}{$h_m$}
    \includegraphics[width=5cm]{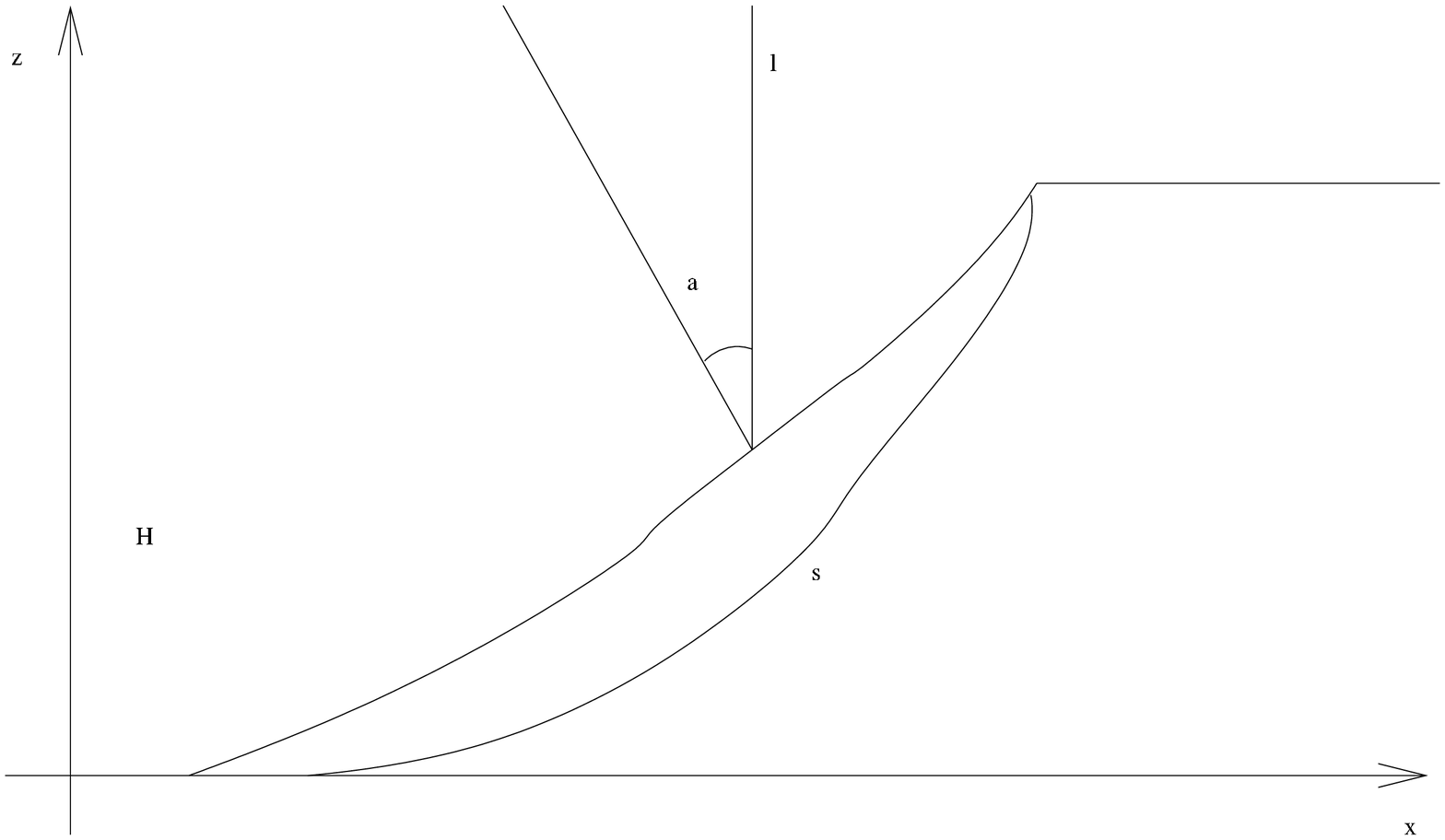}
  \end{psfrags}
  \hfil
  \includegraphics[width=6cm,trim = 0 45 0 0]{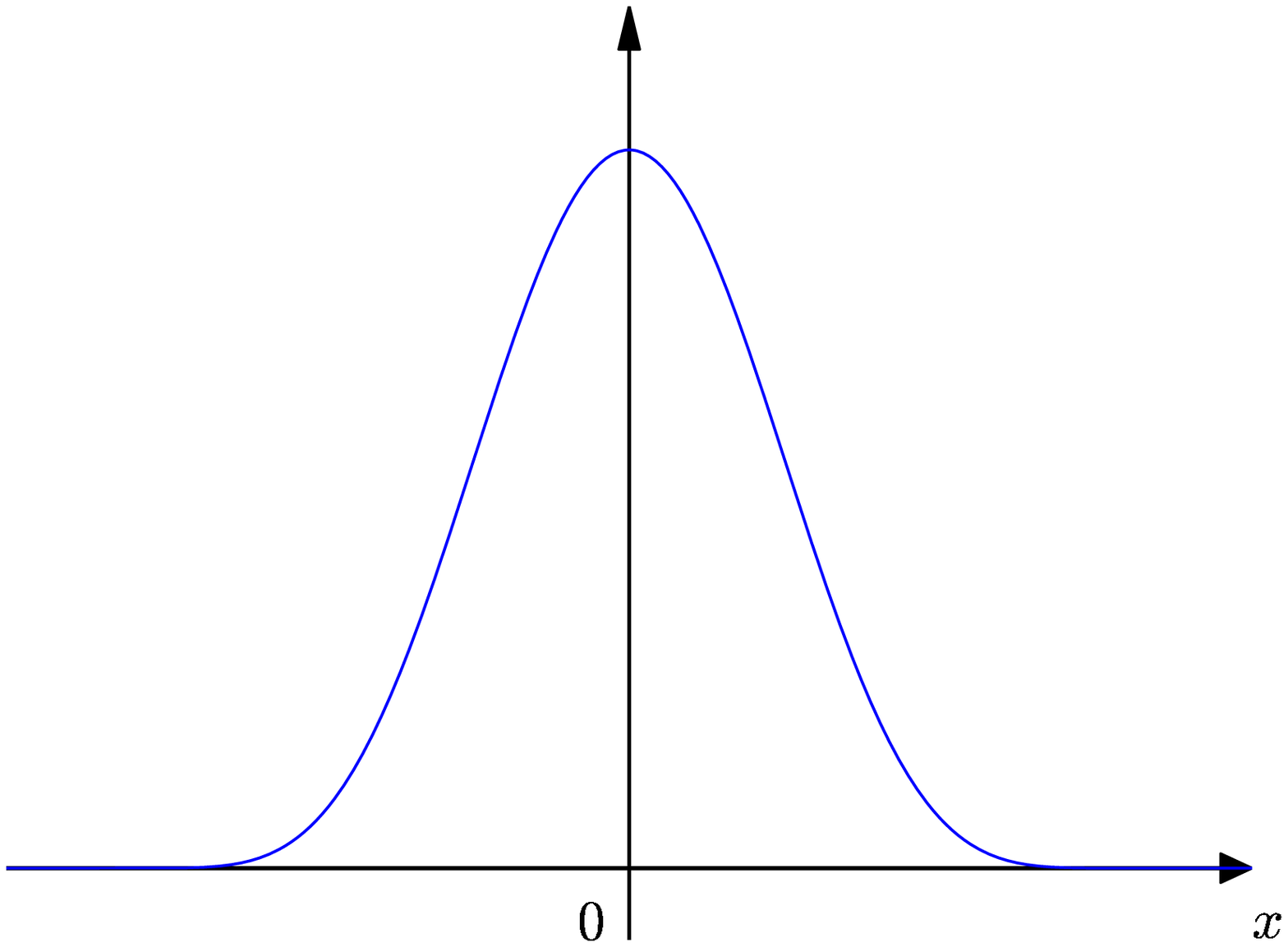}
  \Caption{Left, the incidence angle $\alpha$ of the laser beam on the
    surface $z = H (t,x)$. Right, a possible profile for the functions
    $\mathcal{W}$ and $\mathcal{I}$
    in~\eqref{eq:WL}.\label{fig:angle}}
\end{figure}

A $1$D system of balance laws is used to describe the dynamics of the
melted and of the solid material
in~\cite{ColomboGuerraHertyMarcellini}. Here, we present a description
of this dynamics by means a $2$D system of balance laws of the form:
\begin{equation}
  \label{eq:3}
  \left\{
    \begin{array}{l}
      \partial_t h_m + \div (h_m \, V ) = \mathcal{L}
      \\
      \partial_t h_s = - \mathcal{L}  \,.
    \end{array}
  \right.
\end{equation}
The vector $V = V(t,x)$ describes the projection of the melted
material velocity on the horizontal $(x,y)$-plane. Its modulus must
depend on the wind speed $w = w(t,x)$, which is centered at the laser
beam sited at $x = x_L (t)$. Its direction depends on the geometry of
the melted metal and of the solid surface $z = H (t,x)$, where $H =
h_s + h_m$. The source term $\mathcal{L}$ is directly related to the
laser position and intensity: it describes the net rate at which the
solid part turns into melted. Also $\mathcal{L}$ depends on the metal
geometry, since the heat absorption is strictly related to the
incidence angle between the moving melted metal surface and the
vertical laser beam, see Figure~\ref{fig:angle}, left.

Here, we posit the following assumptions:
\begin{eqnarray}
  \label{eq:V}
  V & = &
  \left(w (t,x) - \tau_g h_m\right) \;
  \frac{-\gradx(\eta*H)}{\sqrt{1+\norma{\gradx(\eta*H)}^2}}
  \\
  \label{eq:I}
  \mathcal{L} & = &
  \frac{i (t,x)}{1+\norma{\gradx(\eta*H)}^2} \,.
\end{eqnarray}
The term with the coefficient $\tau_g$ in~\eqref{eq:V} is related to
the shear stress, inspired by~\cite{ColomboGuerraHertyMarcellini,
  VossenSchuttler}. The denominator in~\eqref{eq:V} is due to a
(smooth) normalization of the direction $-\gradx (\eta * H)$ of the
average steepest descent along the surface $z = H (t,x)$. Indeed, the
convolution kernel $\eta$ is chosen smooth, compactly supported and
with total mass $1$, so that $\gradx\left(\eta * H (t)\right) (x)$ is
the average gradient at position $x$ and time $t$ of the surface $z =
H (t,x)$.

In~\eqref{eq:I}, the numerator $i = i (t,x)$ is related to the laser
intensity. It can be reasonably described through a compactly
supported bell shaped function centered at the location of the moving
focus of the laser beam.
The denominator is the squared cosine of an averaged incidence angle
of the laser on the surface $z = H (t,x)$, see Figure~\ref{fig:angle},
left. In fact,
\begin{displaymath}
  \cos^2 \alpha
  =
  \left(  \frac{[ -\partial_{x_1} H \quad -\partial_{x_2}H \quad 1]
      \left[\begin{array}{c}
          0 \\ 0 \\ 1
        \end{array}\right]}{\norma{[ -\partial_{x_1} H \quad -\partial_{x_2}H \quad 1]} \; \norma{[0\quad 0 \quad 1]}}\right)^2
  =
  \frac{1}{1+\norma{\gradx H}^2} \,.
\end{displaymath}
For the wind function $w = w (t,x)$ and for the laser intensity
function $i = i (t,x)$ we choose a dependence on the form
\begin{equation}
  \label{eq:WL}
  w (t,x) = \mathcal{W}\left(\norma{x-x_L (t)}\right)
  \quad \mbox{ and } \quad
  i (t,x) = \mathcal{I}\left(\norma{x-x_L (t)}\right)
\end{equation}
where both maps $\mathcal{W}$ and $\mathcal{I}$ have the form in
Figure~\ref{fig:angle}, right. More precisely, in the real setting
under consideration, the diameter of the support of $\mathcal{W}$ is a
few times larger than that of $\mathcal{I}$.

\smallskip

We stress that the present model describes how the laser beam digs a
block of metal along its movement, i.e., it describes the dynamics of
the melted metal and the profile of the solid material during the
passing of the laser beam. At the physical level, the actual formation
of the hole makes the melted material fall and, essentially,
disappear. At the analytic level, the appearance of the hole causes
major discontinuities that can hardly be described within a model of
the form~\eqref{eq:3}. Therefore, we
provide~\eqref{eq:3}--\eqref{eq:V}--\eqref{eq:I} with an initial datum
\begin{equation}
  \label{eq:4}
  h_s (0,x) = h_s^o \quad \mbox{ and } \quad h_m (x) = 0
\end{equation}
where the constant $h_m^o$ is the uniform thickness of the plate under
consideration. Then, we interpret the region where $h_s (t,x) < 0$ as
the region where the cut is accomplished.

\smallskip

As a result we obtain the following model:
\begin{equation}
  \label{eq:1}
  \left\{
    \begin{array}{l@{}}
      \displaystyle
      \partial_t h_m
      +
      \div\left[
        \left(w (t,x)h_m - \tau_g (h_m)^2\right)
        \frac{-\gradx(\eta*H)}{\sqrt{1+\norma{\gradx(\eta*H)}^2}}
      \right]
      =
      \frac{i (t,x)}{1+\norma{\gradx(\eta*H)}^2}
      \\
      \displaystyle
      \partial_t h_s
      =
      - \frac{i (t,x)}{1+\norma{\gradx(\eta*H)}^2}
      \\
      H = h_s+h_m
    \end{array}
  \right.
\end{equation}

To apply Theorem~\ref{thm:main} to the model~\eqref{eq:1}, a formal
modification is necessary. Indeed, we introduce a cutoff function
\begin{equation}
  \label{eq:S}
  {\cal T}_g (t,x) = \tau_g \, \mathcal{S}\left(\norma{x-x_L (t)}\right)
  \quad \mbox{ where } \quad
  \mathcal{S} (\xi) =
  \left\{
    \begin{array}{lrcl}
      1 & \xi & \in & [0, \, r]
      \\
      0 & \xi & \in & \left[R, +\infty\right[
    \end{array}
  \right.
\end{equation}
for a smooth $\mathcal{S}$ and suitable (\emph{large}) $r$ and $R$,
with $r <
R$.  % From the modeling point of view, its role is justified
% by the fact that no metals moves far from the laser beam.
We thus obtain
\begin{equation}
  \label{eq:1bis}
  \left\{
    \begin{array}{l@{}}
      \displaystyle
      \partial_t h_m
      +
      \div
      \left[
        \frac{
          -\left(
            w (t,x)h_m - {\cal T}_g (t,x) (h_m)^2
          \right)\gradx(\eta*H)}{\sqrt{1+\norma{\gradx(\eta*H)}^2}}
      \right]
      =
      \frac{i (t,x)}{1+\norma{\gradx(\eta*H)}^2}
      \\
      \displaystyle
      \partial_t h_s
      =
      - \frac{i (t,x)}{1+\norma{\gradx(\eta*H)}^2}
      \\
      H = h_s+h_m
    \end{array}
  \right.
\end{equation}
When used with real data, the two problems~\eqref{eq:1}
and~\eqref{eq:1bis} are indistinguishable.

\begin{proposition}
  \label{prop:OK}
  The model~\eqref{eq:1bis} fits into~\eqref{eq:2} setting:
  \begin{displaymath}
    \begin{array}{@{}r@{\;}c@{\;}l}
      N & = & 2
      \\
      n & = & 2
      \\
      m & = & 2
      \\
      u_1 & = & h_m
      \\
      u_2 & = & h_s
    \end{array}
    \qquad
    \theta (x)
    =
    \left[
      \begin{array}{cc}
        \partial_{x_1} \eta (x) & \partial_{x_1} \eta (x)
        \\
        \partial_{x_2} \eta (x) & \partial_{x_2} \eta (x)
      \end{array}
    \right]
    \qquad
    \begin{array}{@{}r@{\;}c@{\;}l@{}}
      \phi_1 (t,x,u_1,A)
      & = &
      \!\left(w (t,x) - {\cal T}_g (t,x) u_1\right)
      \frac{- u_1 \, A}{\sqrt{1+\norma{A}^2}}
      \\
      \phi_2 (t,x,u_2,A) & = & 0
      \\
      \Phi_1 (t,x,u,A) & = & \frac{1}{\sqrt{1+\norma{A}^2}} \, i (t,x)
      \\
      \Phi_2 (t,x,u,A) & = & - \frac{1}{\sqrt{1+\norma{A}^2}} \, i (t,x) \,,
    \end{array}
  \end{displaymath}
  where $w,i$ are defined in~\eqref{eq:WL} and ${\cal T}_g$
  in~\eqref{eq:S}.  Moreover, if
  \begin{equation}
    \label{eq:uffa}
    x_L \in (\C2 \cap \W{2}{\infty}) ([0,\hat T]; \reali^2)
    \,,\quad
    \mathcal{W}, \mathcal{I}, \mathcal{S} \in \Cc2 (\reali; \reali)
    \quad \mbox{ and } \quad
    \eta \in \Cc3 (\reali^2; \reali)
  \end{equation}
  for a positive $\hat T$, then,
  assumptions~\textbf{($\boldsymbol{\phi}$)},
  \textbf{($\boldsymbol{\Phi}$)}, \textbf{($\boldsymbol{\theta}$)}
  hold.
\end{proposition}

\noindent The proof is deferred to Section~\ref{sec:TD}.  The above
Proposition~\ref{prop:OK} allows to apply Theorem~\ref{thm:main} to
model~\eqref{eq:1bis}, ensuring its well posedness.

\subsection{Numerical Integration}
\label{subs:NI}

The model~\eqref{eq:1}, fed with realistic values of the various
parameters, is able to reproduce the rising of ripples. The following
numerical integrations show this qualitative feature.

We use below the numerical method presented
in~\cite{AggarwalColomboGoatin}, where it is proved to be convergent
up to a subsequence in the case of a system of nonlocal conservation
laws. As it is usual, we deal with the source terms by means of the
fractional step method, see for instance~\cite[Section
12.1]{LeVeque}. In other words, we use a Lax--Friedrichs type
algorithm for the convective part and a first order explicit forward
Euler method for the ordinary differential equations arising from the
source terms.

The computational domain is the rectangle $[0, \, 40] \times [-2, \,
2]$, entirely contained in the metal plate to be cut (all lengths
being measured in millimeters). The mesh size is $5 \cdot 10^{-3}$
along both axis. The integration is computed for $t \in [0, \, 1]$,
time being measured in seconds.

The laser trajectory is
\begin{displaymath}
  x_L (t) =
  \left\{
    \begin{array}{lrcl}
      (3, \, 0) & t & \in [0, \, 0.1]
      \\
      \left(3 + 40 \, (t-0.1), \, 0\right) & t & \in \left]0.1, \, 1\right]
    \end{array}
  \right.
\end{displaymath}
meaning that for $t \in [0, \, 0.1]$ the initial hole is drilled
centered at $(3, \, 0)$, in the interior of the metal plate. The speed
of the laser beam, $40 \frac{{\rm mm}}{{\rm sec}}$, is coherent with
the data in~\cite{lasermelt}, see
also~\cite[Table~1]{ColomboGuerraHertyMarcellini}.
\begin{figure}[!h]
  \centering
  \includegraphics[width=0.33\textwidth, trim = 10 20 25
  30]{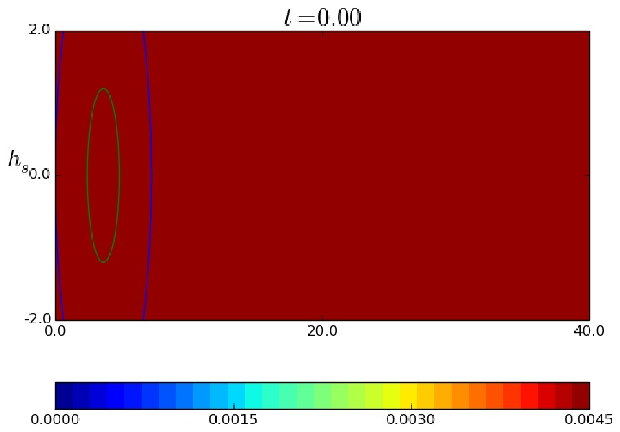}%
  \includegraphics[width=0.33\textwidth, trim = 10 20 25
  30]{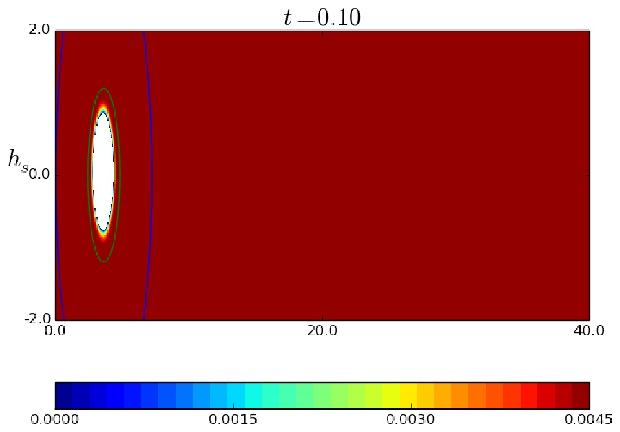}%
  \includegraphics[width=0.33\textwidth, trim = 10 20 25 30]{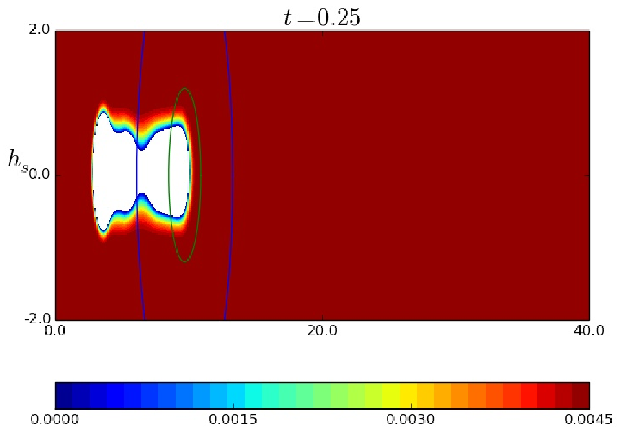}\\
  \includegraphics[width=0.33\textwidth, trim = 10 20 25
  30]{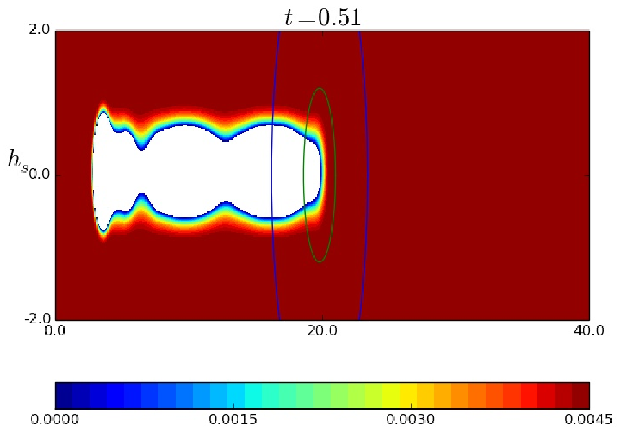}%
  \includegraphics[width=0.33\textwidth, trim = 10 20 25
  30]{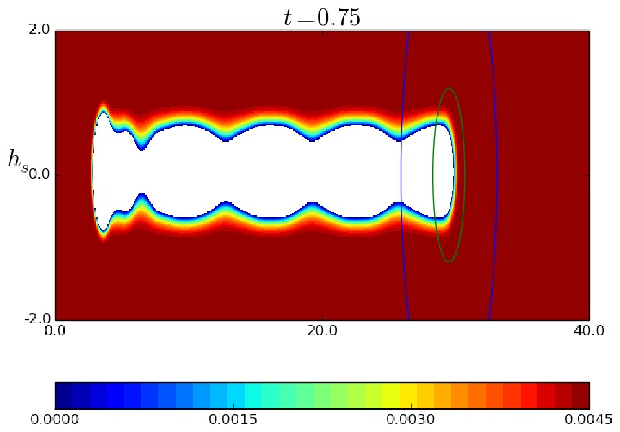}%
  \includegraphics[width=0.33\textwidth, trim = 10 20 25 30]{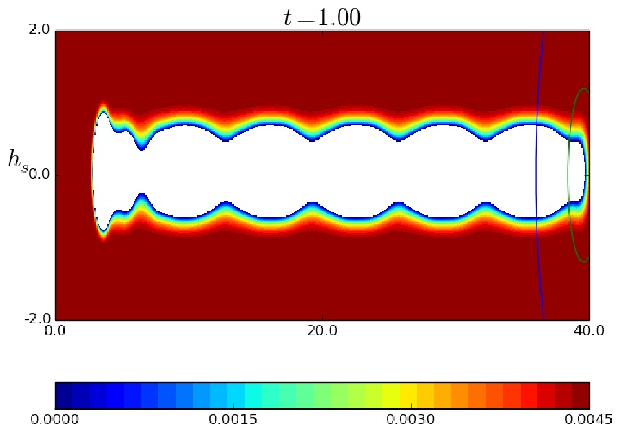}\\
  \caption{Numerical integration of~\eqref{eq:1} with the data and
    parameters provided in~\S~\ref{subs:NI},
    see~\cite{ColomboGuerraHertyMarcellini, lasermelt}. These are the
    contour plots of the solid metal level $h_s$ in the interval $[0,
    \, 4.5]$ over the domain $[0,\, 40] \times [-2, \, 2]$
    millimeters. The inner circle (appearing as an ellipse due to the
    different scales on the two axis) is the support of the laser
    beam. The outer one is the support of the wind. At time $t=0.1$,
    the initial hole is terminated and the laser beam starts moving
    rightwards. Note the formation of \emph{``ripples''}, i.e., the
    sides of the cut are not flat but present an apparently regularly
    oscillating profile. Neither data nor parameters are
    \emph{``pulsating''}.}
  \label{fig:NI}
\end{figure}

The wind and laser functions are given by~\eqref{eq:WL} setting
\begin{displaymath}
  \begin{array}{rclr@{\;}c@{\;}lrcl}
    \mathcal{W} (\xi)
    & = &
    \displaystyle
    \left(1-\left(\frac{\xi}{3.6}\right)^2\right)^4
    &
    \mbox{ for } \norma{\xi} & \leq & 3.6
    \\
    \mathcal{I} (\xi)
    & = &
    \displaystyle
    2 \, \left(1-\left(\frac{\xi}{1.2}\right)^2\right)^6
    &
    \mbox{ for } \norma{\xi} & \leq & 1.2
  \end{array}
\end{displaymath}
corresponding to a laser beam with radius $1.2\,$mm,
see~\cite{lasermelt}. The radius of the surface where the wind blows
downward is $3$ times that of the laser beam. Besides, we set $\tau_g
= 4$. The convolution kernel is
\begin{displaymath}
  \eta (x) =
  \frac{1}{\int_{\reali^2} \tilde \eta (y)\d{y}} \, \tilde \eta (x)
  \quad \mbox{ where } \quad
  \tilde\eta (x) = \left(1-\left(\frac{\norma{x}}{2.4}\right)^2\right)^3
  \qquad
  \mbox{ for } \norma{x} \leq 2.4 \,.
\end{displaymath}
\begin{figure}[!h]
  \centering
  \includegraphics[width=0.33\textwidth]{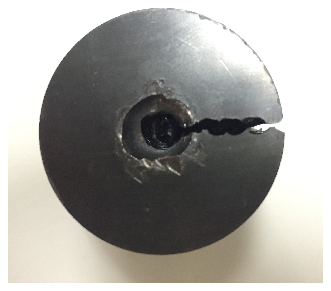}%
  \includegraphics[width=0.67\textwidth]{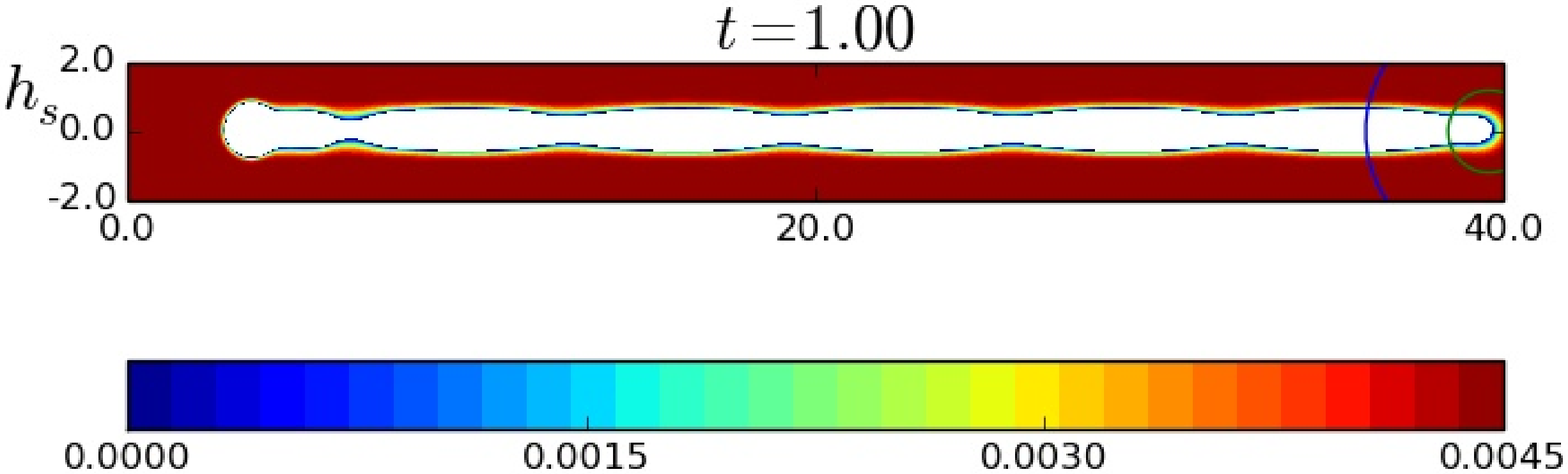}%
  \Caption{Left, a real piece of metal with a hole and a cut made by a
    laser beam. Right, the result of the numerical integration of the
    model~\eqref{eq:1} as in Figure~\ref{fig:NI}, but plotted with the
    same scales along the two axis.}
  \label{fig:real}
\end{figure}
As initial datum, we choose
\begin{displaymath}
  \bar h_m (x) = 0
  \qquad \mbox{ and } \qquad
  \bar h_s (x) = 4.5
  \qquad \mbox{ for all } x \in \reali^2\,,
\end{displaymath}
representing a flat metal plate $4.5\,$mm thick, see~\cite{lasermelt}.

The result of this integration is in Figure~\ref{fig:NI}, which
displays the contour plot of the surface $z = h_s (t,x)$ remaining
after the cut, restricted to the interval $[0, \, 4.5]$mm. The white
part corresponds to a level below $0$ and should be understood as
expelled, corresponding to the cut. Remark the oscillations arisen
along the sides of the cut. No parameter and no datum in the
integration oscillates, nevertheless, the solution displays these sort
of \emph{``ripples''}.

\section{Materials Flowing on a Conveyor Belt}
\label{sec:CB}

\begin{figure}[!h]
  \centering
  \begin{psfrags}
    \psfrag{ml}{$-\ell$} \psfrag{l}{$\ell$} \psfrag{x1}{$x_1$}
    \psfrag{x2}{$x_2$}
    \includegraphics[width=0.45\textwidth]{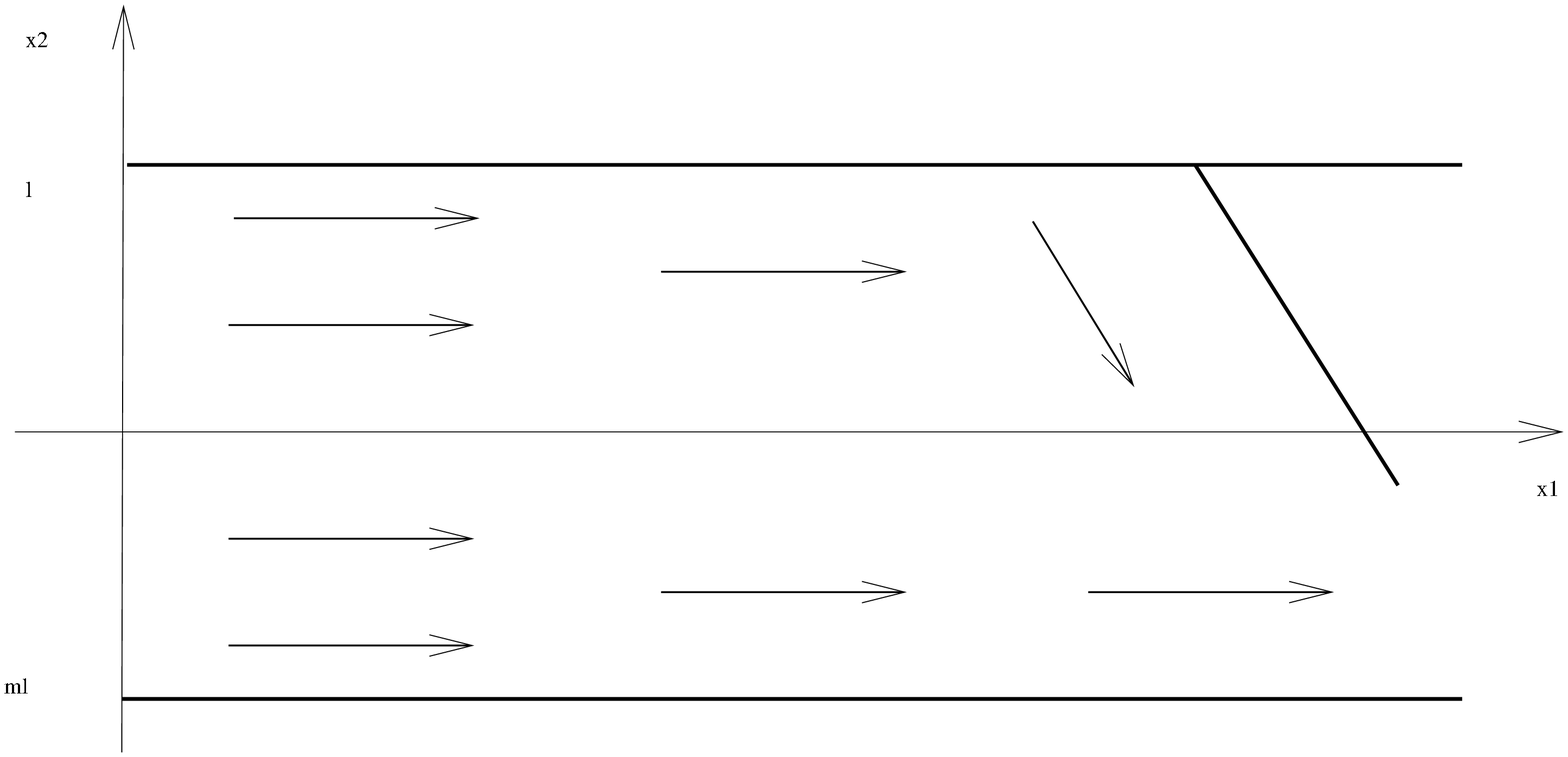}
  \end{psfrags}
  \hfil
  \begin{psfrags}
    \psfrag{ml}{$-\ell$} \psfrag{l}{$\ell$} \psfrag{x1}{$x_1$}
    \psfrag{x2}{$x_2$} \psfrag{d}{$\delta$}
    \includegraphics[width=0.45\textwidth]{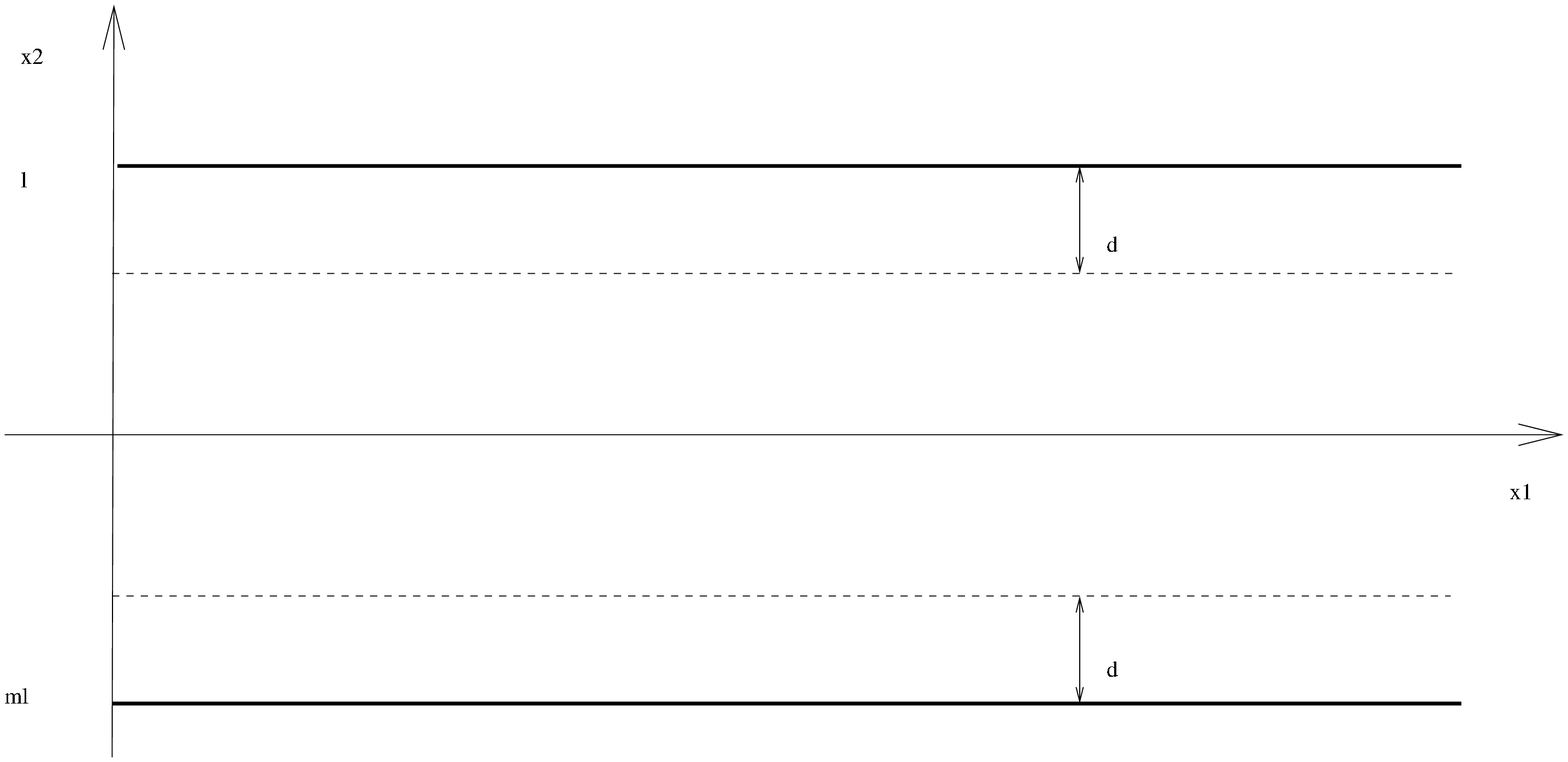}
  \end{psfrags}
  \Caption{Left, a conveyor belt with a selector restricting the
    possible path of the carried cargo. Right, geometry and notation
    of the conveyor belt.}
  \label{fig:conBelt}
\end{figure}

A macroscopic model for the flow of materials along a conveyor belt is
presented in~\cite{SchleperGoettlich}. The material consists of a
large number of solid identical particles, called \emph{cargo}. From a
macroscopic point of view, the cargo state is identified by a density
$\rho = \rho (t,x)$, where $t$ is time and $x \equiv (x_1, \, x_2)$ is
the coordinate along the conveyor belt. The industrial interest behind
these modeling efforts is motivated by the need of an efficient
management of specific parts of the production process. A standard
example is the pouring of newly produced bolts in boxes. In this case,
a \emph{selector} is positioned on the belt to drive the bolts in a
short segment of the belt, so that at the end of the conveyor they
fall in their boxes, see~\cite{SchleperGoettlich} and
figure~\ref{fig:conBelt}, left. For other references on these modeling
issues, both from the microscopic and macroscopic points of view, see
for instance~\cite{Minkin}, related to conveyor belts in mines,
or~\cite{Veronika1, Sekler} and the review~\cite{Jahangirian20101}.

With the notation in~\cite[Section~3]{SchleperGoettlich}, a
macroscopic description for the cargo dynamics is provided by the
equation
\begin{equation}
  \label{eq:conv}
  \partial_t \rho
  +
  \div \left(
    \rho
    \left(
      \boldsymbol{v^{stat}} (x)
      +
      H(\rho - \rho_{max}) \,
      \mathcal{I} (\rho)
    \right)
  \right) = 0 \,.
\end{equation}
Here, $\boldsymbol{v^{stat}}$ is the time independent velocity of the
underlying conveyor belt. The fixed positive $\rho_{max}$ is the
maximal cargo density and $H$ is the usual Heaviside function. The
term
\begin{equation}
  \label{eq:Icargo}
  \mathcal{I} (\rho)
  =
  \epsilon \,
  \frac{-\gradx (\eta * \rho)}{\sqrt{1 + \norma{\gradx (\eta * \rho)}^2}}
\end{equation}
describes how the cargo velocity is modified when the maximal density
is reached: particles move towards regions with lower average cargo
density, $\eta$ being a $\Cc2$ positive function with integral $1$, so
that $\eta * \rho$ is an average cargo density. Further details are
available in~\cite[Section~3]{SchleperGoettlich},
where~\eqref{eq:conv} is supplied with suitable boundary conditions
along the sides of the conveyor belt. The numerical study therein
shows a good agreement between the solutions to~\eqref{eq:conv} and
real data.

Next, we slightly modify~\eqref{eq:conv}. The conveyor belt is
described by the strip $\modulo{x_2} \leq \ell$. First, we replace the
Heaviside function by a regularization
\begin{equation}
  \label{eq:H}
  H^\mu \in \C2 (\reali; [0,1])
  \quad \mbox{ with } \quad
  H^\mu (\xi) = H (\xi) \quad \forall \xi
  \mbox{ with }\modulo\xi > \mu \,.
\end{equation}
Then, we modify $\boldsymbol{v^{stat}} (x)$ so that it incorporates
the upper and lower conveyor boundaries. To this aim, we introduce the
vector field $\boldsymbol{b} (x) \in \Cc2 (\reali^2; \reali^2)$, see
Figure~\ref{fig:conBelt}, right,
such that:
\begin{displaymath}
  \begin{array}{r@{\;}c@{\;}r@{\qquad}r@{\;}c@{\;}l}
    \left(\boldsymbol{b} (x_1, \ell)\right)_2 & = & -\hat \epsilon
    \quad \forall x_1 \in \reali \,,
    &
    \left(\boldsymbol{b} (x)\right)_1 & = & 0
    \quad \forall x \in \reali^2 \,,
    \\
    \left(\boldsymbol{b} (x_1, -\ell)\right)_2 & = & \hat \epsilon
    \quad \forall x_1 \in \reali \,,
    &
    \left(\boldsymbol{b} (x)\right)_2 & = & 0 \quad \forall x \in \reali^2 \mbox{ with } \modulo{x_2-\ell} > \delta \mbox{ or } \modulo{x_2+\ell} > \delta\,.
  \end{array}
\end{displaymath}
We therefore obtain the equation
\begin{equation}
  \label{eq:cargo}
  \partial_t \rho
  +
  \div \left(
    \rho
    \left(
      \boldsymbol{v^{stat}} (x)
      +
      \boldsymbol{b} (x)
      +
      H^\mu(\rho - \rho_{max}) \,
      \mathcal{I} (\rho)
    \right)
  \right) = 0
\end{equation}
which describes the cargo dynamics along the conveyor belt. Thanks to
the framework provided by Theorem~\ref{thm:main}, we can incorporate
in the model also the cargo source and sink. Indeed, we assume that
the solid particles are poured on the belt in a region, say, $R_{in} =
[0, a] \times [-\ell, \ell]$ and fall out of the belt in the region
$R_{out} = [L-a, L] \times [-\ell, \ell]$. To this aim, for a positive
$\hat T$, we introduce the source and sink functions
\begin{equation}
  \label{eq:Psi}
  \begin{array}{@{}r@{\;}c@{\;}l@{\mbox{ with ~}}l@{}}
    \Psi_{in} & \in & \C2 ([0,\hat T] \times R_{in}; \reali^+)
    &
    \spt \Psi_{in} (t, \cdot ) \subseteq R_{in}
    \mbox{ for all } t \in [0,\hat T]
    \\[6pt]
    \Psi{_{out}} & \in & \C2 ([0,\hat T] \times R_{in} \times\reali ; \reali^+)
    &
    \begin{array}{@{}l@{}}
      \spt \Psi_{out} (t, \cdot, \rho ) \subseteq R_{out}
      \mbox{ for all } t \in [0,\hat T]\mbox{ and }\rho \in \reali
      \\
      \Psi_{out} (\cdot , \cdot ,\rho) = 0
      \mbox{ for all } \rho \leq 0
    \end{array}
  \end{array}
\end{equation}
The function $\Psi_{in}$ is the rate at which particles are poured in
$R_{in}$, while $\Psi{_{out}}$ describes the outflow from the belt.

We can assume that the belt is initially empty, thus we obtain the
following Cauchy Problem, where we set $\boldsymbol{v} =
\boldsymbol{v}^{stat} + \boldsymbol{b}$,
\begin{equation}
  \label{eq:PbCargo}
  \left\{
    \begin{array}{l}
      \partial_t \rho
      +
      \div \left(
        \rho
        \left(
          \boldsymbol{v} (x)
          +
          H^\mu(\rho - \rho_{max}) \;
          \mathcal{I} (\rho)
        \right)
      \right)
      =
      \Psi_{in} (t,x) - \Psi_{out} (t,x,\rho)
      \\
      \rho (t,0) =0
    \end{array}
  \right.
\end{equation}

\begin{proposition}
  \label{prop:cargo}
  Fix positive $\hat T, \ell, L, \mu, \rho_{max}, \epsilon, \hat
  \epsilon$ with $\hat \epsilon > \epsilon$.  Let $\mathcal{B} = [0,L]
  \times [-\ell, \ell]$ be the conveyor belt.  If $\boldsymbol{v} \in
  \C2 (\reali^2; \reali^2)$ is such that
  \begin{equation}
    \label{eq:vCargo}
    \begin{array}{rcl@{\qquad}rcl}
      \left(\boldsymbol{v} (0, x_2)\right)_1 & \geq & 0
      \quad \forall x_2 \in [-\ell, \ell]
      &
      \left(\boldsymbol{v} (L, x_2)\right)_1 & \leq & 0
      \quad \forall x_2 \in [-\ell, \ell]
      \\
      \left(\boldsymbol{v} (x_1, -\ell)\right)_2 & \geq & \hat\epsilon
      \quad \forall x_1 \in [0, L]
      &
      \left(\boldsymbol{v} (x_1, \ell)\right)_2 & \geq & \hat\epsilon
      \quad \forall x_1 \in [0, L] \,,
    \end{array}
  \end{equation}
  $H^\mu$ is as in~\eqref{eq:H}, $\mathcal{I}$ is as
  in~\eqref{eq:Icargo} and $\Psi_{in}, \Psi_{out}$ are as
  in~\eqref{eq:Psi}, then there exists a positive $T_*$ such that
  problem~\eqref{eq:PbCargo} admits a solution on the time interval
  $[0,T_*]$. Moreover, this solution is supported in $\mathcal{B}$ for
  all $t \in [0,T_*]$.
\end{proposition}

\noindent The proof is deferred at the end Section~\ref{sec:TD}.

\section{Technical Details}
\label{sec:TD}

The proof of Theorem~\ref{thm:main} consists of several steps. We
briefly describe here the overall formal structure.

Fix a positive $T \in \hat I$ and let $I = [0,T]$. Introduce the map
\begin{displaymath}
  \mathcal{T} \colon (w, \tilde u) \to u
\end{displaymath}
where $u \equiv (u_1, \ldots, u_n)$ and its $i$-th component $u_i$
solves the nonlinear balance law
\begin{equation}
  \label{eq:5}
  \left\{
    \begin{array}{l}
      \partial_t u_i
      +
      \div  \phi_i (t,x,u_i, \theta * w)
      = \Phi_i (t, x, u_i, \theta*w)
      \\
      u_i (0,x) = \tilde u_i (x)
    \end{array}
  \right.
\end{equation}
for $i=1, \ldots, n$.  By construction, solving~\eqref{eq:2} is
equivalent to solving the fixed point problem $u = \mathcal{T} (u,\bar
u)$. The core of the proof thus consists in choosing $T$ and suitable
subsets
\begin{displaymath}
  \mathbb{W} \subset \C0 \left(I;\L1 (\reali^N;\reali^n)\right)
  \quad \mbox{ and } \quad
  \mathbb{U} \subset \L1 (\reali^N;\reali^n),
\end{displaymath}
see~\eqref{eq:defT}, so that
\begin{displaymath}
  \left\{
    \begin{array}{ll@{\qquad}l}
      \mbox{(i)} &
      \forall (w, \tilde u) \in
      \mathbb{W} \times \mathbb{U},
      &
      \mathcal{T} (w,\tilde u)
      \mbox{ is well defined,}
      \\
      \mbox{(ii)} &
      \forall (w, \tilde u) \in
      \mathbb{W} \times \mathbb{U},
      &
      \mathcal{T} (w,\tilde u)
      \mbox{ is in } \mathbb{W} ,
      \\
      \mbox{(iii)} &
      \forall \tilde u \in \mathbb{U},
      &
      w \to \mathcal{T} (w,\tilde u)
      \mbox{ is a contraction,}
      \\
      \mbox{(iv)} &
      \forall w \in \mathbb{W},
      &
      \tilde u \to \mathcal{T} (w,\tilde u)
      \mbox{ is Lipschitz continuous,}
      \\
      \mbox{(v)} &
      \forall (w, \tilde u) \in
      \mathbb{W} \times \mathbb{U},
      &
      t \to \left(\mathcal{T} (w,\tilde u)\right) (t)
      \mbox{ is continuous.}
    \end{array}
  \right.
\end{displaymath}
Steps~\textbf{1} and~\textbf{2} in the proof below give~(i). The
\emph{a priori} bounds proved in steps~\textbf{3}, \textbf{4},
\textbf{5} and~\textbf{6} ensure~(ii). The key estimate~\eqref{eq:f1},
which has the form
\begin{displaymath}
  \norma{\mathcal{T} (w',\tilde u) - \mathcal{T} (w'',\tilde u)}_{\C0 (I;\L1 (\reali^N; \reali^n))}
  \leq \O \; T \; \norma{w'-w''}_{\C0 (I;\L1 (\reali^N; \reali^n))}
\end{displaymath}
and is proved in Step~\textbf{7}, shows that~(iii) holds for $T$
small. The statement~(iv) is obtained in Step~\textbf{9} through an
estimate of the form
\begin{displaymath}
  \norma{\mathcal{T} (w, \tilde u') - \mathcal{T} (w, \tilde u'')}_{\C0 (I;\L1 (\reali^N; \reali^n))}
  \leq \O \; \norma{\tilde u' - \tilde u''}_{\L1 (\reali^N; \reali^n)} \,,
\end{displaymath}
see~\eqref{eq:f2}. Finally, (v) is the content of Step~\textbf{4},
see~\eqref{eq:daqui}--\eqref{eq:aqui}, used also in the proof of~(ii).

Once the statements~(i), $\ldots$, (v) are obtained, the proof of
Theorem~\ref{thm:main} is essentially completed.

\begin{proofof}{Theorem~\ref{thm:main}}
  Throughout, we use the standard properties of the convolution
  product and, in particular, the following bounds. If $\theta$
  satisfies~\textbf{($\boldsymbol{\theta}$)} and $u \in \L1 (\reali^N;
  \reali^n)$, then
  \begin{displaymath}
    \norma{\theta * u}_{\L\infty (I\times \reali^N;\reali^n)}
    \leq
    \norma{\theta}_{\L\infty (\reali^N; \reali^{m\times n})} \;
    \norma{u}_{\C0 (I; \L1 (\reali^N; \reali))}
  \end{displaymath}
  which is a straightforward generalization, for instance,
  of~\cite[Theorem~IV.15]{Brezis}. By~\textbf{($\boldsymbol{\theta}$)},
  without any loss of generality, we may assume that
  \begin{displaymath}
    \norma{\theta_{ji}}_{\L1(\reali^N; \reali)} \leq 1/n
    \quad \mbox{ for all } \quad
    j=1, \ldots, m
    \quad \mbox{ and } \quad
    i = 1, \ldots, n \,.
  \end{displaymath}
  This requirement simplifies several estimates below, since it
  ensures that
  \begin{displaymath}
    u_i (x) \in \mathcal{U}_U \mbox{ for all } i=1, \ldots, n
    \quad \mbox{ and } \quad
    x \in \reali^{N}
    \quad \Rightarrow \quad
    (\theta * u)(x) \in \mathcal{U}_U^m
    \quad \mbox{ for all } \quad x  \in \reali^{N} \,.
  \end{displaymath}

  \paragraph{1: Notation and Definition of $\mathcal{T}$.}

  Fix positive $K$, $U$, $\bar U$, $R$ and $\bar R$ with
  \begin{displaymath}
    \norma{\bar u}_{\L1 (\reali^N; \reali^n)} \leq \bar R < R
    \,,\quad
    \norma{\bar u}_{\L\infty (\reali^N; \reali^n)} \leq \bar U < U
    \quad \mbox{ and } \quad
    \tv (\bar u) < K \,.
  \end{displaymath}
  Introduce the $\L1$ closed sphere centered at the initial datum
  $\bar u$ with radius $R$ and its intersection with $\BV$ as follows:
  \begin{eqnarray*}
    B_{\L1} (\bar u,R,U)
    & = &
    \left\{
      u \in \L1 (\reali^N; \reali^n)
      \colon
      \norma{u - \bar u}_{\L1 (\reali^N; \reali^n)} \leq R
      \mbox{ and }
      u (x) \in \mathcal{U}_U^n
    \right\}
    \\
    B_{\L1 \cap \BV} (\bar u,\bar R,\bar U,K)
    & = &
    \left\{
      u \in B_{\L1} (\bar u, \bar R, \bar U)
      \colon
      \tv (u) \leq K
    \right\} .
  \end{eqnarray*}
  For any positive $T \in \hat I$, denote $I = [0,T]$ and define the
  map
  \begin{equation}
    \label{eq:defT}
    \begin{array}{@{}c@{\;}c@{\;}ccccc@{}}
      \mathcal{T} & \colon &
      \C0\left(I; B_{\L1} (\bar u, R, U)\right)
      & \times &
      B_{\L1\cap \BV} (\bar u, \bar R, \bar U, K)
      & \to
      & \C0\left(I; B_{\L1} (\bar u, R, U)\right)
      \\
      & &
      w &,&\tilde u& \to & u
    \end{array}
  \end{equation}
  where the function $u \equiv (u_1, \ldots, u_n)$ is such that for
  $i=1, \ldots, n$, $u_i$ solves~\eqref{eq:5}. We equip the Banach
  space $\C0\left(I; \L1 (\reali^N; \reali^n)\right)$ with its natural
  norm
  \begin{displaymath}
    \norma{u}_{\C0\left(I; \L1 (\reali^N; \reali^n)\right)}
    =
    \sup_{t \in I} \norma{u (t)}_{\L1 (\reali^N; \reali^n)} \,,
  \end{displaymath}
  and the metric space $B_{\L1\cap \BV} (\bar u, r, K)$ with the
  $\L1$--distance.  Denote below
  \begin{displaymath}
    \Omega_T = I \times \reali^N \times \reali
    \quad \mbox{ and } \quad
    \Omega_T^U = I \times \reali^N \times \mathcal{U}_U \,.
  \end{displaymath}
  Moreover, we set
  \begin{equation}
    \label{eq:Lambda}
    \Lambda (t,U)
    =
    \norma{\lambda (\cdot, \cdot, U)}_{\L1 ([0,t]\times\reali^N; \reali)}
  \end{equation}
  so that $\Lambda (\cdot, U) \in \C0 (\hat I; \reali)$ is non
  decreasing, bounded and $\Lambda (0,U) = 0$ for all $U \in
  \reali^+$.

  Throughout, we denote by $C$ a quantity dependent only on $\lambda$
  and on the norms in~\textbf{($\boldsymbol\phi$)},
  \textbf{($\boldsymbol\Phi$)} and~\textbf{($\boldsymbol\theta$)}, but
  independent of $T$, $R$, $U$, $\bar R$, $\bar U$ and $K$. Similarly,
  $C_U$ is a constant depending only on $\norma{\phi}_{\W2\infty
    (I\times \reali^N \times \mathcal{U}_U\times \mathcal{U}_U^m;
    \reali^{n\times m})}$ and on $\norma{\Phi}_{\W1\infty (I\times
    \reali^N \times \mathcal{U}_U\times \mathcal{U}_U^m;
    \reali^{n})}$.

  \paragraph{2: Problem~\eqref{eq:5} Admits a Solution.}

  Note that if $i \neq j$, equation~\eqref{eq:5} is decoupled from the
  analogous equation for $u_j$. Therefore, we want to apply the
  classical result by Kru\v zkov~\cite[Theorem~1]{Kruzkov}, see
  also~\cite[Theorem~2.1]{MagaliV2}, to each equation in~\eqref{eq:5},
  setting iteratively for $i=1, \ldots, n$
  \begin{displaymath}
    f (t,x,u) = \phi_i \left(t, x, u, (\theta* w) (t,x)\right)
    \quad \mbox{ and } \quad
    F (t,x,u) = \Phi_i \left(t, x, u, (\theta* w) (t,x)\right) \,.
  \end{displaymath}
  To this aim, we check that the assumption~\textbf{(H1*)}
  in~\cite[Theorem~2.1]{MagaliV2}, see also~\cite[Theorem~1]{Kruzkov},
  is satisfied.
  \begin{description}
  \item[(H1*)] $f \in \C0 (\Omega_T; \reali^N)$ holds
    by~\textbf{($\boldsymbol\phi$)}
    and~\textbf{($\boldsymbol\theta$)}, since $w \in \C0\left(I; \L1
      (\reali^N; \reali^n)\right)$.

    $F \in \C0 (\Omega_T; \reali)$ holds
    by~\textbf{($\boldsymbol\Phi$)}
    and~\textbf{($\boldsymbol\theta$)}, since $w \in \C0\left(I; \L1
      (\reali^N; \reali^n)\right)$.

    $f$ has continuous derivatives $\partial_u f$, $\partial_u \gradx
    f$, $\gradx^2 f$, by~\textbf{($\boldsymbol\phi$)}
    and~\textbf{($\boldsymbol\theta$)}.

    $F$ has continuous derivatives $\partial_u F$ and $\gradx F$
    by~\textbf{($\boldsymbol\Phi$)}
    and~\textbf{($\boldsymbol\theta$)}.

    $\partial_u f \in \L\infty (\Omega_U^T; \reali)$
    by~\textbf{($\boldsymbol{\phi}$)}.

    $(F - \div f) \in \L\infty (\Omega_T^U;\reali)$
    by~\textbf{($\boldsymbol\phi$)} and~\textbf{($\boldsymbol\Phi$)}.

    $\partial_u (F - \div f) \in \L\infty (\Omega_T^U;\reali)$
    by~\textbf{($\boldsymbol\phi$)} and~\textbf{($\boldsymbol\Phi$)}.
  \end{description}

  \noindent Therefore, problem~\eqref{eq:5} admits a solution $u \in
  \L\infty \left(I; \Lloc1 (\reali^N; \reali^n)\right)$.

  \paragraph{3: Total Variation Estimate.}

  We want to apply~\cite[Theorem~2.5]{ColomboMercierRosini} as refined
  in~\cite[Theorem~2.2]{MagaliV2}. To this aim, we
  verify~\textbf{(H2*)} in~\cite[\S~2]{MagaliV2}.
  \begin{description}
  \item[(H2*)]$\gradx \partial_u f \in \L\infty
    (\Omega_T^U;\reali^{N\times N})$ by~\textbf{($\boldsymbol\phi$)}
    and~\textbf{($\boldsymbol\theta$)}.

    $\partial_u F \in \L\infty (\Omega_T^U;\reali)$ since $\partial_u
    F = \partial_{u_i}\phi_i$ and since~\textbf{($\boldsymbol\Phi$)}
    holds.

    $\int_I \int_{\reali^N} \norma{\gradx(F-\div f)
      (t,x,\cdot)}_{\L\infty (\mathcal{U}_U;\reali^N)} \d{x} \d{t} <
    +\infty$: indeed, note that the inequality $\norma{\gradx F
      (t,x,\cdot)}_{\L\infty (\mathcal{U}_U;\reali^N)} \leq \lambda
    (t,x,U)$ holds by~\textbf{($\boldsymbol\Phi$)}. Moreover,
    \begin{equation}
      \label{eq:incubo}
      \begin{array}{rcl}
        \div f (t,x,u)
        & = &
        \div \phi_i \left(t,x,u, (\theta * w) (t,x)\right)
        \\
        & &
        +
        \gradA \phi_i \left(t,x,u, (\theta * w) (t,x)\right) \;
        \div (\theta * w) (t,x)
      \end{array}
    \end{equation}
    and passing to the gradient
    \begin{eqnarray*}
      & &
      \gradx \div f (t,x,u)
      \\
      & = &
      \gradx \div \phi_i \left(t,x,u, (\theta * w) (t,x)\right)
      \\
      & &
      +
      \gradA \div \phi_i \left(t,x,u, (\theta * w) (t,x)\right) \;
      \gradx (\theta * w) (t,x)
      \\
      & &
      +
      \gradx \gradA \phi_i \left(t,x,u, (\theta * w) (t,x)\right) \;
      \div (\theta * w) (t,x)
      \\
      & &
      +
      \gradA^2 \phi_i \left(t,x,u, (\theta * w) (t,x)\right) \;
      \gradx (\theta * w) (t,x) \;
      \div (\theta * w) (t,x)
      \\
      & &
      +
      \gradA \phi_i \left(t,x,u, (\theta * w) (t,x)\right) \;
      \gradx \div (\theta * w) (t,x) \,,
    \end{eqnarray*}
    so that, using the standard properties of the convolution
    and~\textbf{($\boldsymbol\phi$)}
    \begin{eqnarray*}
      \norma{\gradx \div f (t,x,u)}
      % & \leq &
      % \Bigl[
      % 1
      % +
      % 3\norma{\gradx \theta}_{\L\infty (\reali^N; \reali^{m\times
      % n})}
      % \norma{w}_{\C0 (I; \L1 (\reali^N; \reali))}
      % \\
      % & &
      % \qquad
      % +
      % \left(\norma{\gradx \theta}_{\L\infty (\reali^N;
      %   \reali^{m\times n})}
      %   \norma{w}_{\C0 (I; \L1 (\reali^N; \reali))}\right)^2
      % \Bigr] \lambda (t,x)
      % \\
      & \leq &
      (1 + 3 \, C_U \, R \, T + C_U \, R^2 \, T^2) \; \lambda (t,x,U)
      \\
      & \leq &
      C_U \;(1 + R\, T + R^2 \, T^2) \; \lambda (t,x,U)
    \end{eqnarray*}
    and hence, using~\textbf{($\boldsymbol{\phi}$)},
    \textbf{($\boldsymbol{\Phi}$)} and~\eqref{eq:Lambda},
    \begin{eqnarray}
      \nonumber
      \!\!\!
      & &
      \int_I \int_{\reali^N} \norma{\gradx(F-\div f)
        (t,x,\cdot)}_{\L\infty (\mathcal{U}_U;\reali^N)} \d{x} \d{t}
      \\
      \nonumber
      \!\!\!
      & \leq&
      \int_I \int_{\reali^N}
      \left(
        \norma{\gradx F (t,x,\cdot)}_{\L\infty (\mathcal{U}_U;\reali^N)}
        +
        \norma{\gradx\div f (t,x,\cdot)}_{\L\infty (\mathcal{U}_U;\reali^N)}
      \right)
      \d{x} \d{t}
      \\
      \nonumber
      \!\!\!
      & \leq &
      \int_I \int_{\reali^N}
      \left(
        \lambda (t,x,U)
        +
        C_U \;(1 + R\, T + R^2 \, T^2) \, \lambda (t,x,U)
      \right)
      \d{x} \d{t}
      \\
      \label{eq:i}
      \!\!\!
      & = &
      C_U \;(1 + R\, T + R^2 \, T^2)\;
      \Lambda (T,U) \,.
    \end{eqnarray}
  \end{description}
  \noindent To apply~\cite[Theorem~2.5]{MagaliV2}, with reference
  to~\cite[(2.6)]{MagaliV2} compute first
  \begin{eqnarray*}
    \gradx \partial_u f (t,x,u)
    & = &
    \gradx \partial_u \phi_i\left(t,x,u,(\theta * w) (t,x)\right)
    \\
    & &
    +
    \gradA \partial_u \phi_i\left(t,x,u,(\theta * w) (t,x)\right)
    \gradx (\theta * w) (t,x)
  \end{eqnarray*}
  so that
  \begin{displaymath}
    \norma{\gradx \partial_u f}_{\L\infty
      (I\times\reali^N\times\mathcal{U}_U;\reali^{N\times N})}
    \leq C_U+ C_U \, C \, R \, T \leq C_U \, (1+C\, R \, T)
  \end{displaymath}
  and
  \begin{eqnarray}
    \nonumber
    \kappa_0^*
    & = &
    (2N+1)
    \norma{\gradx \partial_u f}_{\L\infty
      (I\times\reali^N\times\mathcal{U}_U;\reali^{N\times N})}
    +
    \norma{\partial_u F}_{\L\infty(I\times\reali^N\times\mathcal{U}_U;\reali)}
    \\
    \nonumber
    & \leq &
    (2N+1)\, C_U \, (1 + C \, R\, T) + C_U
    \\
    \label{eq:k0}
    & \leq &
    C \; C_U \; (1+R\,T)\,.
  \end{eqnarray}
  Denoting $W_N = \int_0^{\pi/2} (\cos\theta)^N \d\theta$,
  use~\eqref{eq:i}, \eqref{eq:k0} to obtain, for all $t \in I$,
  \begin{eqnarray}
    \nonumber
    \tv \left(u_i(t)\right)
    & \leq &
    \tv(\tilde u) e^{\kappa_0^* t}
    +
    N \, W_N \!
    \int_I \! e^{\kappa_0^* (t-\tau)} \!\!
    \int_{\reali^N}
    \norma{\gradx(F-\div f) (\tau,x,\cdot)}_{\L\infty (\mathcal{U}_U;\reali^N)}
    \d{x} \d{\tau}
    \\
    \nonumber
    & \leq &
    \left(
      \tv(\tilde u)
      +
      N \, W_N
      \int_I
      \int_{\reali^N}
      \norma{\gradx(F-\div f) (\tau,x,\cdot)}_{\L\infty (\mathcal{U}_U;\reali^N)}
      \d{x} \d{\tau}
    \right)
    e^{\kappa_0^* t}
    \\
    \label{eq:quella}
    & \leq &
    \left(
      K
      +
      C \, C_U \, (1 + R \, T + R^2 \, T^2) \;
      \Lambda (T,U)
    \right)
    e^{C C_U (1+RT) T}  \,.
  \end{eqnarray}

  \paragraph{4: $\L1$ Continuity in Time}

  We use~\cite[Corollary~2.4]{MagaliV2}. To this aim, verify first
  that $f,F$ satisfy~\textbf{(H3*)}.
  \begin{description}
  \item[(H3*)] $\partial_u f \in \L\infty (\Omega^U_T;\reali^N)$,
    already verified in~\textbf{(H1*)}.

    $\partial_u F \in \L\infty (\Omega^U_T;\reali)$, already verified
    in~\textbf{(H2*)}.

    $\int_I \int_{\reali^N} \norma{(F - \div f) (t,x,\cdot)}_{\L\infty
      (\mathcal{U}_U;\reali)} \d{x} \d{t} < +\infty$:
    use~\eqref{eq:incubo}, \textbf{($\boldsymbol{\phi}$)}
    and~\textbf{($\boldsymbol{\Phi}$)} to obtain the bound
    \begin{eqnarray}
      \nonumber
      & &
      \int_I \int_{\reali^N} \norma{(F - \div f)
        (t,x,\cdot)}_{\L\infty (\mathcal{U}_U;\reali)} \d{x} \d{t}
      \\
      \nonumber
      & \leq &
      \int_I \int_{\reali^N}
      \norma{F (t,x,\cdot)}_{\L\infty (\mathcal{U}_U;\reali)} \d{x} \d{t}
      +
      \int_I \int_{\reali^N}
      \norma{\div f (t,x,\cdot)}_{\L\infty (\mathcal{U}_U;\reali)} \d{x} \d{t}
      \\
      \nonumber
      & \leq &
      \Lambda (T,U)
      +
      \int_I \int_{\reali^N}
      \norma{\div \phi \left(t,x,\cdot,(\theta* w) (t,x)\right)}_{\L\infty (\mathcal{U}_U;\reali)}
      \d{x} \d{t}
      \\
      \nonumber
      & &
      +
      \int_I \int_{\reali^N}
      \norma{\gradA \phi \left(t,x,\cdot,(\theta* w) (t,x)\right)}_{\L\infty (\mathcal{U}_U;\reali)}
      \d{x} \d{t}
      \\
      \nonumber
      & &
      \qquad
      \times
      \norma{\div \theta}_{\L\infty (I\times\reali^N;\reali^n)}
      \norma{w}_{\C0 (I;\L1 (\reali^N;\reali))}
      \\
      \nonumber
      & \leq &
      \Lambda (T,U)
      +
      \Lambda (T,U)
      \left(
        1
        +
        \norma{\div \theta}_{\L\infty (I\times\reali^N;\reali^n)}
        \norma{w}_{\C0 (I;\L1 (\reali^N;\reali))}
      \right)
      \\
      \label{eq:questa}
      & = &
      C \, (1 + R \, T) \, \Lambda (T,U)\,.
    \end{eqnarray}
  \end{description}
  \noindent Repeating the same computations on the time interval
  between $s$ and $t$, by~\textbf{($\boldsymbol{\phi}$)},
  \cite[(2.8)]{MagaliV2}, \eqref{eq:quella} and~\eqref{eq:questa}, for
  all $t,s \in I$,
  \begin{eqnarray}
    \label{eq:daqui}
    & &
    \norma{u_i (t) - u_i (s)}_{\L1 (\reali^N;\reali)}
    \\
    \nonumber
    & \leq &
    \!\!
    \modulo{\int_s^t \int_{\reali^N}
      \norma{(F-\div f) (\tau,x,\cdot)}_{\L\infty (\mathcal{U}_U;\reali)}
      \d{x}}
    +
    \modulo{t-s} \;
    \norma{\partial_u f}_{\L\infty (I\times\reali^N\times\mathcal{U}_U,\reali)} \;
    \sup_{\tau \in I} \tv \left(u_i (\tau)\right)
    \\
    \nonumber
    & \leq &
    \!\!
    C (1 + R T) \modulo{\Lambda (t,U) - \Lambda (s,U)}
    \\
    \label{eq:aqui}
    & &
    +
    C \modulo{t-s} \!
    \left[
      K
      +
      C (1 + R T + R^2 T^2) \Lambda (T,U)
    \right] \!
    e^{C C_U (1+RT) T}
  \end{eqnarray}
  proving the uniform $\L1$--continuity in time of the map $t \to u_i
  (t)$, where $u = \mathcal{T} (w, \tilde u)$.

  \paragraph{5: $\L\infty$ Bound}

  Passing to the limit $\epsilon \to 0$ in the classical
  estimate~\cite[Formula~(4.6)]{Kruzkov}, we have that,
  using~\cite[Formul\ae~(4.1), (4.2) and 4)~in~\S~4]{Kruzkov},
  $\modulo{u (t,x)} \leq (M_o + c_o T) e^{c_1 T}$, where
  \begin{eqnarray*}
    M_o
    & = &
    \norma{\bar u (x)}_{\L\infty{(\reali^N; \reali)}}
    \\
    & \leq &
    \bar R \,.
    \\
    c_o
    & = &
    \norma{\div f (\cdot,\cdot,0) - F (\cdot, \cdot, 0)}_{\L\infty (I\times\reali^N;\reali)}
    \\
    & \leq &
    \norma{\div \phi_i \left(\cdot, \cdot, 0, (\theta * w) (\cdot, \cdot)\right)}_{\L\infty (I\times\reali^N;\reali)}
    \\
    & &
    +
    \norma{\gradA \phi_i \left(\cdot, \cdot, 0, (\theta * w) (\cdot, \cdot)\right) \;
      \div (\theta * w) (\cdot, \cdot)}_{\L\infty (I\times\reali^N;\reali)}
    \\
    & &
    +
    \norma{\Phi \left(\cdot, \cdot, 0, (\theta * w) (\cdot, \cdot)\right)}_{\L\infty (I\times\reali^N;\reali)}
    \\
    & \leq &
    C_U + C_U\, C \, R \, T + C_U
    \\
    & \leq &
    C \, C_U \, (1 + R \, T) \,.
    \\
    c_1
    & = &
    \sup_{I\times \reali^N \times \mathcal{U}_U}
    \left(
      - \partial_u \div f (t, x, u)
      +
      \partial_u F (t, x, u)
    \right)
    \\
    & \leq &
    \norma{\partial_u \div \phi_i \left(\cdot, \cdot, \cdot, (\theta * w) (\cdot, \cdot)\right)}_{\L\infty (I\times\reali^N\times\mathcal{U}_U;\reali)}
    \\
    & &
    +
    \norma{\partial_u \gradA \phi_i \left(\cdot, \cdot, \cdot, (\theta * w) (\cdot, \cdot)\right) \;
      \div (\theta * w) (\cdot, \cdot)}_{\L\infty (I\times\reali^N\times\mathcal{U}_U;\reali)}
    \\
    & &
    +
    \norma{\partial_u \Phi \left(\cdot, \cdot, \cdot, (\theta * w) (\cdot, \cdot)\right)}_{\L\infty (I\times\reali^N\times\mathcal{U}_U;\reali)}
    \\
    & \leq &
    C_U + C_U\, C \, R \, T + C_U
    \\
    & \leq &
    C \, C_U \, (1 + R \, T) \,.
  \end{eqnarray*}
  Therefore,
  \begin{equation}
    \label{eq:infty}
    \norma{u}_{\L\infty (I\times\reali^N;\reali)}
    \leq
    \left(
      \bar R + C \, C_U \, (1 + R \, T) T
    \right)
    \exp \left( C \, C_U \, (1 + R \, T) T\right) \,.
  \end{equation}

  \paragraph{6: $\mathcal{T}$ is Well Defined}

  Apply~\eqref{eq:daqui}--\eqref{eq:aqui} with $s=0$, obtaining that
  for all $t \in I$
  \begin{eqnarray*}
    \norma{u_i (t) - \bar u_i}_{\L1 (\reali^N;\reali)}
    & \leq &
    \norma{u_i (t) - \tilde u_i}_{\L1 (\reali^N;\reali)}
    +
    \norma{\tilde u_i - \bar u_i}_{\L1 (\reali^N;\reali)}
    \\
    & \leq &
    C (1 + R T) \Lambda (T)
    +
    C T
    \left[
      K
      +
      C (1 + R T + R^2 T^2) \Lambda (T)
    \right]
    e^{C (1+R T) T}
    +
    \bar R\,.
  \end{eqnarray*}
  This inequality, together with~\eqref{eq:infty}, ensures that if $T$
  is sufficiently small, $u (t) = \left(\mathcal{T} (w,\tilde
    u)\right) (t) \in B_{\L1} (\bar u, R)$ for all $t \in I$. This
  estimate, together with what was proved at~\textbf{2}
  and~\textbf{4}, ensures that $\mathcal{T}(w,\tilde u) \in \C0
  \left(I;\L1 (\reali^N; \reali^n)\right)$ for any $\tilde u \in
  B_{\L1\cap \BV} (\bar u,r,K)$.

  \paragraph{7: $\mathcal{T}$ is a Contraction.}

  Here we use the stability
  result~\cite[Theorem~2.6]{ColomboGuerraHertyMarcellini} as refined
  in~\cite[Theorem~2.5]{MagaliV2}. To this aim, for $w',w'' \in
  \C0\left(I; \L1 (\reali^N; \reali^n)\right)$, call $f',f'',F',F''$
  the corresponding fluxes and sources. We first verify that $f'-f''$
  and $F'-F''$ satisfy~\textbf{(H3*)}.
  \begin{description}
  \item[(H3*)] $\partial_u (f'-f'') \in \L\infty
    (\Omega_T^U;\reali^N)$ is proved as in~\textbf{(H1*)}.

    $\partial_u (F'-F'') \in \L\infty (\Omega_T^U;\reali)$ is proved
    as in~\textbf{(H2*)}.

    Using~\textbf{($\boldsymbol{\Phi}$)} and~\eqref{eq:incubo},
    \begin{eqnarray}
      & &
      \int_I \int_{\reali^N}
      \norma{(F'-F'')-\div (f'-f'') (t,x,\cdot)}_{\L\infty
        (\mathcal{U}_U;\reali)} \d{x} \d{t}
      \\
      \nonumber
      & \leq &
      \int_I \int_{\reali^N}
      \norma{(F'-F'') (t,x,\cdot)}_{\L\infty
        (\mathcal{U}_U;\reali)} \d{x} \d{t}
      \\
      \nonumber
      & + &
      \int_I \int_{\reali^N}
      \norma{\div (f'-f'') (t,x,\cdot)}_{\L\infty
        (\mathcal{U}_U;\reali)} \d{x} \d{t}
      \\
      \nonumber
      & \leq &
      \int_I \int_{\reali^N}
      \norma{
        \Phi_i\left(t,x,\cdot, (\theta*w') (t,x)\right)
        -
        \Phi_i\left(t,x,\cdot, (\theta*w'') (t,x)\right)
      }_{\L\infty (\mathcal{U}_U;\reali)}
      \d{x} \d{t}
      \\
      \nonumber
      & + &
      \int_I \int_{\reali^N}
      \norma{
        \div \left[
          \phi_i \left(t,x,\cdot, (\theta* w') (t,x)\right)
          -
          \phi_i \left(t,x,\cdot, (\theta* w'') (t,x)\right)
        \right]
      }_{\L\infty (\mathcal{U}_U;\reali)}
      \d{x} \d{t}
      \\
      \nonumber
      & + &
      \int_I \int_{\reali^N}
      \Big\Vert
      \gradA \phi_i \left(t,x,\cdot, (\theta* w') (t,x)\right)
      \div (\theta* w') (t,x)
      \\
      \nonumber
      & &
      \qquad\qquad\qquad
      -
      \gradA \phi_i \left(t,x,\cdot, (\theta* w'') (t,x)\right)
      \div (\theta* w'') (t,x)
      \Big\Vert_{\L\infty (\mathcal{U}_U;\reali)}\d{x}\d{t}
      \\
      \nonumber
      & \leq &
      \norma{\Phi}_{\W{1}{\infty} (I \times \reali^N \times \mathcal{U}_U \times \mathcal{U}_U^m; \reali^n)}
      \int_I \int_{\reali^N}
      \norma{\left(\theta* (w'-w'')\right) (t,x)}
      \d{x} \d{t}
      \\
      \nonumber
      & + &
      \norma{\phi}_{\W{2}{\infty} (I \times \reali^N \times \mathcal{U}_U \times \mathcal{U}_U^m; \reali^{n\times N})}
      \int_I \int_{\reali^N}
      \norma{\left(\theta* (w'-w'')\right) (t,x)}
      \d{x} \d{t}
      \\
      \nonumber
      & + &
      \int_I \int_{\reali^N}
      \Big\Vert
      \Big[
      \gradA \phi_i \left(t,x,\cdot, (\theta* w') (t,x)\right)
      \\
      \nonumber
      & &
      \qquad\qquad\qquad
      -
      \gradA \phi_i \left(t,x,\cdot, (\theta* w'') (t,x)\right)
      \Big]
      \div (\theta* w') (t,x)
      \Big\Vert_{\L\infty (\mathcal{U}_U;\reali)}\d{x}\d{t}
      \\
      \nonumber
      & + &
      \int_I \int_{\reali^N}
      \Big\Vert
      \gradA \phi_i \left(t,x,\cdot, (\theta* w'') (t,x)\right)
      \\
      \nonumber
      & &
      \qquad\qquad\qquad
      \left[
        \div (\theta* w') (t,x)
        -
        \div (\theta* w'') (t,x)
      \right]
      \Big\Vert_{\L\infty (\mathcal{U}_U;\reali)}\d{x}\d{t}
      \\
      \nonumber
      & \leq &
      \norma{\Phi}_{\W{1}{\infty} (I \times \reali^N \times \mathcal{U}_U \times \mathcal{U}_U^m; \reali^n)}
      \norma{\theta}_{\L\infty (\reali^N; \reali^{m\times n})}
      \norma{w' - w''}_{\C0 (I;\L1 (\reali^N; \reali^n))}
      T
      \\
      \nonumber
      & + &
      \norma{\phi}_{\W{2}{\infty} (I \times \reali^N \times \mathcal{U}_U \times \mathcal{U}_U^m; \reali^{n\times N})}
      \norma{\theta}_{\L\infty (\reali^N; \reali^{m\times n})}
      \norma{w' - w''}_{\C0 (I;\L1 (\reali^N; \reali^n))}
      T
      \\
      \nonumber
      & + &
      \norma{\phi}_{\W{2}{\infty} (I \times \reali^N \times \mathcal{U}_U\times \mathcal{U}_U^m; \reali^{n\times N})}
      \\
      \nonumber
      & &
      \qquad
      \times
      \int_I \int_{\reali^N}
      \norma{\left(\theta* (w'-w'')\right) (t,x)}
      \d{x} \d{t}
      \norma{\div (\theta * w')}_{\L\infty (I\times\reali^N;\reali^n)}
      \\
      \nonumber
      & + &
      \norma{\phi}_{\W{1}{\infty} (I \times \reali^N \times \mathcal{U}_U \times \mathcal{U}_U^m; \reali^{n\times N})}
      \int_I \int_{\reali^N}
      \norma{\left(\div\theta* (w'-w'')\right) (t,x)}
      \d{x} \d{t}
      \\
      \label{eq:mitica}
      & \leq &
      C \, C_U \, T \, (1+R)
      \norma{w' - w''}_{\C0 (I;\L1 (\reali^N; \reali^n))} \,.
    \end{eqnarray}
  \end{description}
  \noindent Recall the following quantities
  from~\cite[(2.10)]{MagaliV2} and use~\eqref{eq:incubo}:
  \begin{eqnarray}
    \kappa^*
    & = &
    \norma{\partial_u F'}_{\L\infty (I \times\reali^N \times \reali; \reali)}
    +
    \norma{\partial_u \div (f''-f')}_{\L\infty (I \times\reali^N \times \reali; \reali)}
    \\
    \nonumber
    & \leq &
    \norma{\Phi}_{\W{1}{\infty} (\reali^+ \times \reali^N \times \mathcal{U}_U \times \mathcal{U}_U^m; \reali^{n\times N})}
    \\
    \nonumber
    & &
    +
    \norma{\partial_u \div \left(
        \phi_i\left(\cdot,\cdot,\cdot,(\theta*w') (\cdot,\cdot)\right)
        -
        \phi_i\left(\cdot,\cdot,\cdot,(\theta*w'') (\cdot,\cdot)\right)
      \right)}_{\L\infty (I \times\reali^N \times \mathcal{U}_U; \reali)}
    \\
    \nonumber
    & &
    +
    \Big\Vert
    \partial_u
    \gradA \phi_i \left(\cdot,\cdot,\cdot, (\theta * w') (\cdot,\cdot)\right)
    \div (\theta * w') (\cdot,\cdot)
    \\
    \nonumber
    & &
    \qquad\qquad
    -
    \partial_u
    \gradA \phi_i \left(\cdot,\cdot,\cdot, (\theta * w'') (\cdot,\cdot)\right)
    \div (\theta * w'') (\cdot,\cdot)
    \Big\Vert_{\L\infty (I \times\reali^N \times \mathcal{U}_U; \reali)}
    \\
    \nonumber
    & \leq &
    \norma{\Phi}_{\W{1}{\infty} (\reali^+ \times \reali^N \times \mathcal{U}_U \times \mathcal{U}_U^m; \reali^{n\times N})}
    +
    2 \norma{\phi}_{\W{2}{\infty} (\reali^+ \times \reali^N \times \mathcal{U}_U \times \mathcal{U}_U^m; \reali^{n\times N})}
    \\
    \nonumber
    & &
    +
    2
    \norma{\phi}_{\W{2}{\infty} (\reali^+ \times \reali^N \times \mathcal{U}_U \times \mathcal{U}_U^m; \reali^{n\times N})}
    \norma{\div \theta}_{\L\infty (\reali^N;\reali^n)}
    \norma{w'}_{\C0 (I;\L1 (\reali^N; \reali^n))}
    \\
    \nonumber
    & &
    +
    \norma{\phi}_{\W{2}{\infty} (\reali^+ \times \reali^N \times \mathcal{U}_U \times \mathcal{U}_U^m; \reali^{n\times N})}
    \norma{\div \theta}_{\L\infty (\reali^N;\reali^n)}
    \norma{w'-w''}_{\C0 (I;\L1 (\reali^N; \reali^n))}
    \\
    \nonumber
    & \leq &
    3 \, C_U + 4 \, C_U \, C \, R \, T
    \\
    \label{eq:k}
    & \leq &
    C \, C_U \, (1+R \, T) \,.
    \\
    \nonumber
    M
    & = &
    \norma{\partial_u f''}_{\L\infty (I \times\reali^N \times \reali; \reali)}
    \\
    \nonumber
    & = &
    \norma{\partial_u
      \phi_i\left(\cdot,\cdot,\cdot,(\theta*w'') (\cdot,\cdot)\right)
    }_{\L\infty (I \times\reali^N \times \mathcal{U}_U; \reali)}
    \\
    \nonumber
    & \leq &
    \norma{\phi}_{\W{2}{\infty} (\reali^+ \times \reali^N \times \mathcal{U}_U \times \mathcal{U}_U^m; \reali^{n\times N})}
    \\
    \nonumber
    & \leq &
    C_U \,.
  \end{eqnarray}
  By~\cite[Remark~2.8]{MagaliV2}, \eqref{eq:k0} and~\eqref{eq:k}
  \begin{equation}
    \label{eq:frac}
    \frac{e^{\kappa_0^* t} - e^{\kappa^* t}}{\kappa_0^* - \kappa^*}
    \leq
    t \, e^{\max\{\kappa_0^*, \kappa^*\}t}
    \leq
    t \, e^{C C_U (1+RT)t}
  \end{equation}
  so that we can prepare the bound
  \begin{eqnarray}
    \nonumber
    & &
    \norma{\partial_u (f'-f'')}_{\L\infty (I\times \reali^N\times \mathcal{U}_U)}
    \\
    \nonumber
    & = &
    \norma{
      \partial_u \phi_i
      \left(\cdot, \cdot, \cdot, (\theta * w') (\cdot, \cdot)\right)
      -
      \partial_u\phi_i
      \left(\cdot, \cdot, \cdot, (\theta * w'') (\cdot, \cdot)\right)
    }_{\L\infty (I\times \reali^N\times \mathcal{U}_U)}
    \\
    \nonumber
    & \leq &
    \norma{\phi_i}_{\W2\infty (I\times \reali^N \times \mathcal{U}_U \times \mathcal{U}_U^m;\reali^N)} \,
    \norma{\theta}_{\L\infty (\reali^N; \reali^{n\times m})} \,
    \norma{w'-w''}_{\C0 (I;\L1 (\reali^N;\reali))}
    \\
    \label{eq:basta}
    & \leq &
    C \, C_U \, \norma{w'-w''}_{\C0 (I;\L1 (\reali^N;\reali))} \,,
  \end{eqnarray}
  and we can finally pass to the key estimate provided
  by~\cite[Theorem~2.5]{MagaliV2} using~\eqref{eq:frac},
  \eqref{eq:basta}, \eqref{eq:i}, \eqref{eq:k} and~\eqref{eq:mitica}
  \begin{eqnarray}
    \nonumber
    & &
    \norma{u_i' (t) - u_i'' (t)}_{\L1 (\reali^N;\reali)}
    \\
    & \leq &
    \frac{e^{\kappa_0^* t} - e^{\kappa^* t}}{\kappa_0^* - \kappa^*} \,
    \tv (\tilde u) \,
    \norma{\partial_u (f'-f'')}_{\L\infty (I\times \reali^N\times \mathcal{U}_U)}
    \\
    \nonumber
    & &
    +
    N W_N
    \int_0^t \frac{e^{\kappa_0^* (t-\tau)} - e^{\kappa^* (t-\tau)}}{\kappa_0^* - \kappa^*}
    \int_{\reali^N}
    \norma{\gradx (F' - \div f') (\tau,x,\cdot)}_{\L\infty (\mathcal{U}_U;\reali^N)}
    \d{x} \d{\tau}
    \\
    \nonumber
    & &
    \qquad
    \times
    \norma{\partial_u (f'-f'')}_{\L\infty (I\times \reali^N\times \mathcal{U}_U)}
    \\
    \nonumber
    & &
    +
    \int_0^t
    e^{\kappa^* (t-\tau)}
    \int_{\reali^N}
    \norma{\left((F'-F'') - \div (f'-f'')\right) (t,x,\cdot)}_{\L\infty (\mathcal{U}_U;\reali^N)}
    \d{x} \d{\tau}
    \\
    \nonumber
    & \leq &
    t \, e^{C C_U (1+RT)t} \, K \, C \, C_U \,
    \norma{w'-w''}_{\C0 (I;\L1 (\reali^N;\reali))}
    \\
    \nonumber
    & &
    +
    C  \, e^{C C_U (1+RT)t}
    C_U (1 + R\, T + R^2 \, T^2) \, \Lambda (T,U)
    \, C \, C_U \,
    \norma{w'-w''}_{\C0 (I;\L1 (\reali^N;\reali))}
    \\
    \nonumber
    & &
    +
    e^{C C_U (1+RT)}
    C \, C_U \, T \, (1+R) \,
    \norma{w'-w''}_{\C0 (I;\L1 (\reali^N;\reali))}
    \\
    \label{eq:f1}
    & \leq &
    C \, C_U \, T
    \left(
      1
      +
      K
      +
      (1 + R\, T + R^2 \, T^2) \, \Lambda (T,U)
      +
      R
    \right)
    e^{C C_U (1+RT)} \,
    \norma{w'-w''}_{\C0 (I;\L1 (\reali^N;\reali))}
  \end{eqnarray}
  which shows that there exists a positive $T_*$, such that the map
  \begin{displaymath}
    \begin{array}{ccccc}
      \mathcal{T}_{\tilde u} & \colon &
      \C0\left([0,T_*]; B_{\L1} (\bar u, R,U)\right) &
      \to &
      \C0\left([0,T_*]; B_{\L1} (\bar u, R,U)\right)
      \\
      & & w & \to & \mathcal{T} (w, \tilde u)
    \end{array}
  \end{displaymath}
  is a contraction, for any $\tilde u \in B_{\L1\cap \BV} (\bar u,
  \bar R, \bar U, K)$.

  \paragraph{8: The Fixed Point of $\mathcal{T}$ Is the Unique
    Solution to~\eqref{eq:2}.}

  The fixed point of $\mathcal{T}_{\tilde u}$ solves~\eqref{eq:2} by
  Definition~\ref{def:sol} and from~\eqref{eq:defT}. On the other
  hand, any solution to~\eqref{eq:2} in the sense of
  Definition~\ref{def:sol}, is a fixed point of $\mathcal{T}_{\tilde
    u}$, proving also uniqueness.

  \paragraph{9: Continuous Dependence on the Initial Datum.}

  Note first that $\mathcal{T}$ is $\L1$-Lipschitz continuous in its
  second argument. Indeed, applying
  again~\cite[Theorem~2.5]{MagaliV2}, we have:
  \begin{eqnarray}
    \nonumber
    \norma{
      \mathcal{T} (w, \tilde u')
      -
      \mathcal{T} (w,\tilde u'')}_{\C0 (I; \L1 (\reali^N; \reali^n))}
    & \leq &
    e^{\kappa^* T} \, \norma{\tilde u' - \tilde u''}_{\L1 (\reali^N; \reali^n)}
    \\
    \label{eq:f2}
    & \leq &
    e^{C C_U (1+R T)}
    \norma{\tilde u' - \tilde u''}_{\L1 (\reali^N; \reali^n)} \,.
  \end{eqnarray}
  By~\cite[Theorem~2.7]{BressanLectureNotes}, the $\L1$--Lipschitz
  continuous dependence on the fixed point of $\mathcal{T}_{\tilde u}$
  from $\tilde u$ follows.
\end{proofof}

\begin{proofof}{Proposition~\ref{prop:OK}}
  To prove that \eqref{eq:1bis} fits into the class~\eqref{eq:2},
  simply observe that
  \begin{equation}
    \label{eq:theta}
    \theta * u = \gradx \left(\eta* (u_1+u_2)\right) \,.
  \end{equation}
  The regularity required in~\textbf{($\boldsymbol{\phi}$)}
  and~\textbf{($\boldsymbol{\Phi}$)} is immediate, the cutoff function
  ${\cal T}_g$ being useful in bounding the terms $\gradA \phi$ and
  $\gradA^2 \phi$. The various estimates follow from~\eqref{eq:uffa}
  and from the fact that the map $(x_1, x_2) \to (x_1,x_2)/ \sqrt{1+
    {x_1}^2 + {x_2}^2}$ is bounded, with all first and second
  derivatives also bounded.
\end{proofof}

\begin{proofof}{Proposition~\ref{prop:cargo}}
  Observe first that~\eqref{eq:vCargo}, the assumption $\hat \epsilon
  > \epsilon$ and~\eqref{eq:Icargo} ensure that the flow in the
  convective part of~\eqref{eq:PbCargo} points inward all along the
  boundary of $\mathcal{B}$. Therefore, if there is a solution
  to~\eqref{eq:PbCargo}, its support is contained in $\mathcal{B}$ for
  all times.  To apply Theorem~\ref{thm:main}, we introduce a function
  $s \in \Cc2 (\reali^2; \reali)$ such that $s (x)=1$ for all $x \in
  \mathcal{B}$. Then, note that~\eqref{eq:PbCargo} belongs to the
  class~\eqref{eq:2}. Indeed, similarly to~\eqref{eq:theta}, set
  \begin{equation}
    \label{eq:OKcargo}
    \begin{array}{@{}r@{\;}c@{\;}l}
      N & = & 2
      \\
      n & = & 1
      \\
      m & = & 2
      \\
      u & = & \rho
    \end{array}
    \quad
    \theta (x)
    =
    \left[
      \begin{array}{@{}cc@{}}
        \partial_{x_1} \eta (x) & \partial_{x_1} \eta (x)
        \\
        \partial_{x_2} \eta (x) & \partial_{x_2} \eta (x)
      \end{array}
    \right]
    \quad
    \begin{array}{@{}r@{\;}c@{\;}l@{}}
      \phi(t,x,u,A) & = &
      u \left(
        \boldsymbol{v} (x)
        -
        \frac{\epsilon H^\mu (\rho-\rho_{max}) \, A}{\sqrt{1+\norma{A}^2}}
      \right)
      s (x)
      \\
      \Phi (t,x,u,A) & = &
      \Psi_{in} (t,x) - \Psi_{out} (t,x,u)\,,
    \end{array}
  \end{equation}
  The invariance of $\mathcal{B}$ proved above ensures that the
  function $s$ has no effect whatsoever on the dynamics described
  by~\eqref{eq:PbCargo}. Therefore, with the given initial datum (as
  well as with any other initial datum supported in $\mathcal{B}$),
  any solution to~\eqref{eq:2}--\eqref{eq:OKcargo} also
  solves~\eqref{eq:PbCargo}, and \emph{viceversa}. The estimates
  required in~\textbf{($\boldsymbol{\phi}$)}
  and~\textbf{($\boldsymbol{\Phi}$)} now immediately follow.
\end{proofof}

\noindent\textbf{Acknowledgment:} Both authors thank M.~Herty (RWTH) and Markus Nie\ss en (Fraunhofer ILT) for several useful discussions. The present work was supported by the PRIN~2012 project \emph{Nonlinear Hyperbolic Partial Differential Equations, Dispersive and Transport Equations: Theoretical and Applicative Aspects} and by the GNAMPA~2014 project
\emph{Conservation Laws in the Modeling of Collective Phenomena}.

{\small

  \bibliographystyle{abbrv}

  \bibliography{laser2_revised}

}

\end{document}